\documentclass[12pt]{article}
\usepackage{latexsym}
\usepackage{amssymb}
\font\gorditas = msbm8
\def\bbb#1{\hbox {{\gordas #1}}}
\def\errita{\hbox{\gorditas R}}

\font\gordas = msbm10 at 12pt
\def\bbb#1{\hbox {{\gordas #1}}}
\def\erre{{\bbb R}}

\def\UNO{1\mkern-7mu1}

\newtheorem{theorem}{Theorem}[section]

\newtheorem{corollary}[theorem]{Corollary}
\newtheorem{proposition}[theorem]{Proposition}

\newtheorem{remark}[theorem]{Remark}

\begin{document}
\noindent

\begin{center}
{\bf\large Occupation times of branching systems\\
 with initial inhomogeneous
 Poisson states \\[.1cm]
and related superprocesses$^*$}\\[1cm]
\begin{tabular}{ccc}
Tomasz Bojdecki$^1$& Luis G. Gorostiza$^\dagger$&Anna  Talarczyk$^1$\\
tobojd@mimuw.edu.pl&lgorosti@math.cinvestav.mx&annatal@mimuw.edu.pl
\end{tabular}
\end{center}
\setcounter{section}{1}
\setcounter{equation}{0}
\footnote{\kern-.6cm $^*$ Supported in part by CONACyT Grant 45684-F (Mexico) and MEiN Grant 1P03A01129 (Poland).\\
$^1$ Institute of Mathematics, University of Warsaw, ul. Banacha 2, 02-097 Warszawa, Poland.\\
$^\dagger$ Centro de Investigaci\'on y de Estudios Avanzados, A.P. 14-740, Mexico 07000 D.F., Mexico.}
\vglue.5cm
\noindent
\centerline{
{\bf Abstract}}
\vglue.5cm
The $(d,\alpha,\beta,\gamma)$-branching particle system consists of particles 
moving in $\erre^d$ according to a symmetric $\alpha$-stable L\'evy process $(0<\alpha\leq 2)$, splitting with a critical  $(1+\beta)$-branching law $(0<\beta\leq 1)$, and starting from an inhomogeneous Poisson random measure with intensity measure 
$\mu_\gamma(dx)=dx/(1+|x|^\gamma), \gamma\geq 0$.  
By means of time rescaling $T$ and Poisson intensity measure $H_T\mu_\gamma$,
occupation time fluctuation limits for the system as $T\to\infty$ have been obtained in two special cases: Lebesgue measure ($\gamma=0$, the homogeneous case),
 and  finite  measures $(\gamma>d)$. 
In some cases  $H_T\equiv 1$ and in others $H_T\to\infty$ as $T\to\infty$ (high density systems).
The limit processes are quite different for Lebesgue and for finite measures. Therefore the question arises of what kinds of limits can be obtained for Poisson intensity measures that are intermediate between Lebesgue measure and finite measures. In this paper the measures $\mu_\gamma, \gamma\in (0,d]$, are used for investigating this question. Occupation time fluctuation limits are obtained which interpolate in some way between the two previous extreme cases. The limit processes depend on different arrangements of the parameters $d,\alpha,\beta,\gamma$. There are two thresholds for the dimension $d$. The first one, $d=\alpha/\beta+\gamma$, determines the need for high density or not in order to obtain non-trivial limits, and its relation with a.s. local extinction of the system is discussed. The second one, $d=[\alpha(2+\beta)-\gamma\vee \alpha)]/\beta$\ (if  $\gamma<d$), interpolates between the two extreme cases, and it is a critical dimension which separates different qualitative behaviors of the limit processes, in particular long-range dependence in ``low'' dimensions, and independent increments in ``high'' dimensions. 
In low dimensions the temporal part of the limit process is a new self-similar stable process which has two different long-range dependence regimes depending on relationships 
among
the  parameters.
 Related results for the corresponding $(d,\alpha,\beta,\gamma)$-superprocess are also given.
\vskip.5cm
\noindent
{\bf AMS 2000 subject classifications:} Primary 60F17, Secondary 60J80, 60G18, 60G52.
\\
{\bf Key words}: Branching particle system, superprocess, occupation time fluctuation, limit theorem, stable process, long-range dependence.
\newpage
\setcounter{section}{1}
\setcounter{equation}{0}
\noindent
{\bf 1 Introduction}
\vglue.5cm
Occupation time fluctuation limits have been proved for the  so-called $(d,\alpha,\beta)$-branching particle systems in $\erre^d$ with initial Poisson 
states in two special cases, namely,
if the Poisson intensity measure is 
either Lebesgue measure, denoted by $\lambda$, or a finite measure \cite{BGT1}, \cite{BGT2}, 
\cite{BGT3}, \cite{BGT4}, \cite{BGT6}. Those cases are 
quite
 special, as explained below, and   the limit processes are 
very
 different. Therefore 
the question arises of what happens with Poisson intensity measures that are intermediate between Lebesgue measure and finite measures. That is the main motivation for the present paper, and our aim is to obtain limit processes that interpolate in some way between those of the two special cases. One of our  objectives
 is to find out when the limits have long-range dependence behavior and to describe it.
Another
motivation is to derive analogous results for the corresponding superprocesses.

In a
$(d,\alpha,\beta)$-branching particle system   the particles move independently in $\erre^d$ according to a standard spherically symmetric $\alpha$-stable L\'evy process, $0<\alpha\leq 2$, the particle lifetime is exponentially distributed with parameter $V$, and the branching law is critical with generating function
\begin{equation}
\label{eq:1.1}
s+\frac{1}{1+\beta}(1-s)^{1+\beta},\quad 0<s<1,
\end{equation}
where $0<\beta\leq 1$ (called $(1+\beta)$-branching law), which is binary branching for $\beta=1$.
The parameter $V$ is not particularly relevant, but it is convenient to use it.
 The empirical measure process $N=(N_t)_{t\geq 0}$ is defined by 
\begin{equation}
\label{eq:1.2}
N_t(A)=\,\hbox{\rm number of particles in the Borel set}\,\, A
\subset\erre^d\,\,\hbox{\rm at time}\,\, t.
\end{equation}
A common assumption for the initial distribution $N_0$ is to take a Poisson random measure with  locally finite intensity measure $\mu$. The corresponding $(d,\alpha,\beta)$-superprocess $Y=(Y_t)_{t\geq 0}$ is a measure-valued process, which is a high-density/short-life/small-particle limit of the particle system, with $Y_0=\mu$. See \cite{D}, \cite{E}, \cite{P2} for background on those particle systems and superprocesses. 
In this paper we investigate (the limiting behavior of)
 the corresponding occupation time processes, i.e.,
$$\int^t_0N_sds,\, t\geq 0,\,\,{\rm and}\,\,\int^t_0Y_sds,\,\, t\geq 0.$$
We recall 
that the distributions of these processes
are characterized by their Laplace functionals as follows \cite{GLM}, \cite{DP2}:
\begin{equation}
\label{eq:1.3}
E{\rm exp}\left\{-\int^t_0\langle N_s,\varphi\rangle ds\right\}={\rm exp}\{-\langle\mu,v_\varphi(t)\rangle\},\,\,\,\varphi\in{\cal S}(\erre^d),
\end{equation}
where $v_\varphi(x,t)$ is the unique (mild) solution of the non-linear equation
\begin{eqnarray}
\label{eq:1.4}
\frac{\partial}{\partial t}v_\varphi&=&\Delta_\alpha v_\varphi-\frac{V}{1+\beta}v^{1+\beta}_\varphi+\varphi(1-v_\varphi),\\
v_\varphi(x,0)&=&0,\nonumber
\end{eqnarray}
and
\begin{equation}
\label{eq:1.5}
E{\rm exp}\left\{-\int^t_0\langle Y_s,\varphi\rangle ds\right\}={\rm exp}\{-\langle\mu,u_\varphi(t)\rangle\},\,\,\,\varphi\in{\cal S}(\erre^d),
\end{equation}
where $u_\varphi(x,t)$ is the unique (mild) solution of the non-linear equation
\begin{eqnarray}
\label{eq:1.6}
\frac{\partial}{\partial t}u_\varphi&=&\Delta_\alpha u_\varphi-\frac{V}{1+\beta}u^{1+\beta}_\varphi+\varphi,\\
u_\varphi(x,0)&=&0,\nonumber
\end{eqnarray}
and $\Delta_\alpha$ is the infinitesimal generator of the $\alpha$-stable process. (See the end of the Introduction for the standard notations $\langle\quad, \quad\rangle,{\cal S}(\erre^d)$.)

For $\mu=\lambda$,  $N_0$ is homogeneous Poisson. This case is special (and technically simpler) because $\lambda$ is invariant for the $\alpha$-stable process (which implies in particular that  $EN_t=\lambda$ for all $t$), and there is the following persistence/extinction dichotomy  \cite{GW}, 
which
heuristically
 explains the need for high density in some cases in order to obtain non-trivial occupation time fluctuation limits, and anticipates the situation we will encounter in this paper:

\noindent
(i) Persistence: If $d>\alpha/\beta$, then $ N_t$ converges in law to an equilibrium state $N_\infty$ as $t\to\infty$, such that $EN_\infty=\lambda$.

\noindent
(ii) Extinction: If $d\leq \alpha/\beta$, then $ N_t$ becomes locally extinct in probability as $t\to\infty$, i.e., for any bounded Borel set $A,N_t(A)\to 0$
in probability. 

\noindent
An analogous persistence/extinction dichotomy holds for the corresponding superprocess (with $Y_0=\lambda$) \cite{DP1}. 
For $\alpha=2$ (Brownian motion) and $d<2/\beta$   a stronger extinction holds: the superprocess becomes locally extinct in finite time a.s. \cite{I2}, and we shall see that so does the particle system.

The case of  $\mu$  finite is special because the particle system  goes to extinction globally in finite time a.s. for every dimension $d$, and so does the superprocess  \cite{P2}.

The time-rescaled occupation time fluctuation process 
$X_T=(X_T(t))_{t\geq 0}$ of the particle system is defined by
\begin{equation}
\label{eq:1.7}
X_T(t)=\frac{1}{F_T}\int^{Tt}_0(N_s-EN_s)ds,
\end{equation}
where $N_s$ is given by (\ref{eq:1.2}) and  $F_T$ is a norming. The problem is to find $F_T$ such that $X_T$ converges in distribution (in some way) as $T\to\infty$, and to identify the limit process and study its properties. This was done for $\mu=\lambda$ in the persistence case, $d>\alpha/\beta$,
\cite{BGT3}, \cite{BGT4}. For the extinction case, $d\leq\alpha/\beta$,
 in \cite{BGT6} we introduced high density, meaning that the initial Poisson intensity measure was taken  of the form $H_T\lambda$, with $H_T\to\infty$ as $T\to\infty$, so as to counteract the local extinction, and we obtained similarly high-density limits with  $\mu$ finite. For $\mu=\lambda$ and $d>\alpha/\beta$ the same results 
hold with or without high density (with different normings) \cite{BGT6}. 
The limit processes are  different for $\mu=\lambda$ and $\mu$ finite (some differences are mentioned below), and the results for any finite measure are essentially the same.

In order to study  asymptotics of $X_T$ as $T\to\infty$ with  measures $\mu$ that are intermediate between the two previous cases, we consider Poisson intensity measures  of the form
\begin{equation}
\label{eq:1.8}
\mu_\gamma(dx)=\frac{1}{1+|x|^\gamma}dx,\quad \gamma\geq 0.
\end{equation}
We call the  model so defined a $(d,\alpha,\beta,\gamma)$- branching particle system, and a $(d,\alpha,\beta,\gamma)$-superpro\-cess the corresponding measure-valued process. 
To obtain non-trivial limits we multiply $\mu_\gamma$ by $H_T$, which is suitably chosen in each case.
Since $\mu_0=\lambda$ and $\mu_\gamma$ is finite for $\gamma> d$, by varying $\gamma$ in the interval $(0,d]$ we  obtain limits  that are between those of the two previous cases, which  are extreme in this sense, in a way that interpolates between them. 
The substantial role of $\gamma$ was already noted in the simpler model 
of particle systems without branching \cite{BGT5}.

The above mentioned behaviors of $N_t$ in the cases $\gamma=0 $ and $\gamma>d$
 raise the following questions on what happens for $\gamma\in (0,d]$, and on its effect on  asymptotics of $X_T$: When does $N_t$ suffer {\it a.s. local extinction} in the sense that for each bounded Borel set $A$ there is a finite random time $\tau_A$ such that $N_t(A)=0$ for all $t\geq \tau_A$ a.s.? In this case the total occupation time $\int^\infty_0N_t(A)dt$ is finite a.s., and therefore high density is needed in order to obtain non-trivial  limits for $X_T$. For $\gamma>0, N_t(A)$ converges to $0$ in probability as $t\to\infty$ for any bounded Borel set $A$ and every dimension $d$, 
so  local extinction in probability occurs,
but the total occupation time may or may not be finite. It turns out that the threshold between the need for high density or not is given by $d=\alpha/\beta+\gamma$, and then a natural question is whether $d=\alpha/\beta+\gamma$ is also the border to a.s. local extinction of the particle system. We will come back to this question.

The  limits for $X_T$  in \cite{BGT6} are of three  different kinds for both $\mu=\lambda$ and $\mu$ finite. In the first case there is a critical dimension, $d_c=\alpha(1+\beta)/\beta$.
For the ``low'' dimensions, $d<d_c$, the limit has a simple spatial structure (the measure $\lambda$) and a complex temporal structure (with long-range dependence). For the ``high'' dimensions, $d>d_c$, the limit has a complex spatial structure (distribution-valued) and a simple temporal structure (with stationary independent increments). For the ``critical'' dimension, $d=d_c$, the spatial and the temporal structures are both simple, but the order of the fluctuations $(F_T)$ is larger, as is typical in phase transitions. The limit processes are always continuous for $d<d_c$, and for $d\geq d_c$ they are continuous if and only if $\beta=1$ (when the limits are Gaussian). For $\mu$ finite, an analogous trichotomy of   results holds, with a new critical dimension, $d_c=\alpha(2+\beta)/(1+\beta)$, 
 another difference being that the limits for the critical and high dimensions are constant in time for $t>0$.

In this paper we show analogous  limits of
 $X_T$ for $(d,\alpha,\beta,\gamma)$-branching particle system;
the critical dimension changes between the ones above, $\alpha(1+\beta)/\beta$ for $\gamma=0$, and $\alpha(2+\beta)/(1+\beta)$ for $\gamma>d$, and they are linked with a unified formula, 
\begin{equation}
\label{eq:1.9}
d_c(\gamma)=\alpha\frac{2+\beta}{\beta}-\frac{\gamma\vee\alpha}{\beta},
\end{equation}
which interpolates between the two  cases (see Remark 2.2(a) for a precise statement). There are several limit processes depending on different arrangements of $d, \alpha, \beta, \gamma$. Some  are analogous to those for $\mu=\lambda$, and some are similar  to those for $\mu$ finite (or even essentially the same). For $\gamma<d$ there are six different cases that include the three ones  recalled above for $\gamma=0$. For $\gamma>d$ there are the three cases obtained in 
\cite{BGT6} (generally for finite $\mu$). In the case $\gamma<d$ and 
$d<d_c(\gamma)$, the temporal structure of the limit is a new real, stable, self-similar, continuous, long-range dependence process $\xi$, defined in (2.1) below, which has two different long-range dependence regimes if $\alpha <2$ (Theorem 1(a) and  Proposition 2.3). This 
strange type of  long-range dependence behavior  already appears in the homogeneous  case, $\gamma=0$ \cite{BGT3}, \cite{BGT6}. An analogous phenomenon occurs with $0<\gamma <d$, the border between the two long-range dependence regimes changes continuously with $\gamma$, and it disappears in the limit $\gamma \nearrow d$ (see formula (2.13)). For $\gamma>d$ there is
only one long-range dependence regime, not depending on $\gamma,\beta$.

In  \cite{BGT1}, \cite{BGT2}, \cite{BGT3}, \cite{BGT4}, \cite{BGT5},
\cite{BGT6}
 we have given the convergence results for $X_T$ in a strong
form (functional convergence when it holds), but in the present article our main objective is identifying the limits,  so
 we have not attempted to prove the strongest form of convergence in each case, 
nevertheless
 we expect that convergence in law in a space of continuous functions holds in all cases where the limit is continuous. We prove functional convergence only in the case of the above mentioned long-range dependence process $\xi$ because of its special properties. A technical difficulty for the tightness proof is the lack of moments if $\beta<1$.

The time-rescaled occupation time fluctuation process for the $(d,\alpha,\beta,\gamma)$-superprocess $Y$ is defined analogously as (\ref{eq:1.7}), 
\begin{equation}
\label{eq:1.10}
X_T(t)=\frac{1}{F_T}\int^{Tt}_0(Y_s-EY_s)ds,
\end{equation}
 and the limits are  obtained  from (the proofs of) the results for the $(d,\alpha,\beta,\gamma)$-branching particle systems, as a consequence of the fact that the log-Laplace equation of the occupation time of the superprocess is simpler than that of the particle system 
(see  (\ref{eq:1.3}), (\ref{eq:1.4}), and (\ref{eq:1.5}), 
(\ref{eq:1.6}).)

Our results
on the fluctuation limits of superprocesses generalize those of Iscoe 
\cite{I1}, who considered the homogeneous case $(\gamma=0)$ only.

Let us come back to the question of high density and local extinction  for the 
$(d,\alpha,\beta,\gamma)$-branching particle system.
From Theorems 2.1, 2.5 and 2.6 it follows immediately (see Corollary~2.10) that in all cases where high density is not necessary (i.e., we may take $H_T\equiv 1$) there is no a.s.\ local extinction (in spite of the fact that local extinction in probability occurs if $\gamma>0$). For instance,
 condition (\ref{eq:2.6})
 in  Theorem 2.1(a) holds automatically if $d>\alpha/\beta+\gamma$ (with $\gamma < d)$, hence high density is not necessary for a non-trivial limit of $X_T$ in this case.  On the other hand, high density is indispensable to obtain a non-trivial limit if either $d\leq \alpha/\beta+\gamma$ or $\alpha<\gamma\leq d$. In the latter case the total occupation time of any bounded Borel set by the process $N$ is finite a.s. (it has finite mean). 
We  prove a.s. local extinction for $\alpha =2$  and  $d<2/\beta+\gamma$, and
we conjecture that a.s. local extinction holds generally for $d<\alpha/\beta+\gamma$ also if $\alpha<2$. This conjecture is supported by the fact that for $d=\alpha/\beta+\gamma, \gamma<\alpha$, there is an ergodic result  (Proposition 2.9). 
For  
$\alpha=2$ and $d<2/\beta +\gamma$ it follows from \cite{I2} (Theorem $3_\beta$) that  the 
$(d,2,\beta,\gamma)$-superprocess  suffers a.s. local extinction. 
The method of \cite{Z} can also be used to prove this  (private communication).
The proof of a.s. local extinction of the particle system consists in showing
that a.s. local extinction of the $(d,2,\beta,\gamma)$-superprocess implies a.s. local extinction of the $(d,2,\beta,\gamma)$-branching particle system (Theorem 2.8). 
 This  implication is not as
 simple as it might seem because the well-known Cox relationship between the particle system and the superprocess (i.e., for each $t$,  $N_t$ is a doubly stochastic Poisson random measure with random intensity
measure given by 
 $Y_t$) is not enough to relate the long time behaviors of the two processes.
But the argument  does not work for $\alpha<2$. In this case the superprocess $Y$ has the instantaneous propagation of support property, i.e., with probability 1 for each $t>0$ if the closed support of $Y_t$ is not empty, then it is all of $\erre^d$.
This follows from the result proved in \cite{P1} for finite initial measure and finite variance branching $(\beta=1$ in our model), which is  extended in 
\cite{LZ} for more general superprocesses and branching mechanisms (including $\beta<1$ in our case). 
It follows that  for the $(d,\alpha,\beta,\gamma)$-superprocess with 
$\alpha<2$, a.s. local extinction and global extinction are equivalent, and it is known that if the initial measure has infinite total mass, the probability of global extinction in finite time is 
$0$. Nevertheless, the total occupation time of a bounded set for the superprocess with $\alpha<2$ may or may not be finite, and this is what is directly relevant for us (see the proof of Theorem 2.8). Iscoe \cite{I1} showed that for 
initial Lebesgue measure and $\alpha=2$ the total occupation time of a bounded set is finite if and only if $d<2/\beta$, and he conjectured that an analogous result holds for $\alpha<2$ and $d<\alpha/\beta$. So far as we know, this conjecture has not been proved.

Summarizing, if  $\gamma<d$ and $\gamma<\alpha $, 
 there are two thresholds for the asymptotics of 
$X_T$, namely, $\alpha/\beta+\gamma$, and $d_c(\gamma)$ given by 
(\ref{eq:1.9}).  The first one, which is smaller than the second one, 
appears to be the border to a.s. local extinction (we know that it is for $\alpha=2$),  and it determines the need for high density.  The second one is the critical dimension between changes of behavior of the limit processes, in particular the change from long-range dependence to independent increments, and from continuity to discontinuity if $\beta<1$.
An interpretation of $d_c(\gamma)$ in terms of the model seems rather mysterious, even in the case $\gamma=0$ (see \cite{BGT4}, Section 4, for several questions on the  meaning of results).

The general methods of proof   developed in 
\cite{BGT3}, \cite{BGT4}, \cite{BGT6}, and the special cases for $\beta=1$, where the limits are Gaussian \cite{BGT1}, \cite{BGT2}, 
can be used for the proofs involving $\mu_\gamma$. 
 However, a considerable amount of  technical work is unavoidable in order to deal with $\gamma>0$. Moreover, each case requires  different calculations.
 We will abbreviate the proofs as much as possible.

Related work appears in \cite{BZ}, \cite{M1}, \cite{M2}
 for the special case $\gamma=0$,  $\beta=1$. \cite{BZ} studies occupation time fluctuations of a single point for a system of binary branching random walks on the lattice with state dependent branching rate. \cite{M1}, \cite{M2} consider general critical finite variance branching laws.
\cite{BZ}, \cite{M1}, \cite {M2} and \cite{M3}  also study  the systems in equilibrium.
We have already
 mentioned the paper \cite{I1} on occupation time fluctuation limits of $(d,\alpha,\beta)$-superprocesses.
 Some other papers regarding extinction, ergodicity and occupation times of branching particle systems and superprocesses are  \cite{BZ}, \cite{CG},
\cite{DGW},  \cite{DR}, \cite{DF},  \cite{EK}, \cite{FG}, \cite{FVW}, \cite{H}, \cite{IL}, \cite{K},
\cite{LR}, \cite{MR}, \cite{M3}, \cite{Sh}, \cite{Ta}, \cite{VW},
\cite{Zh} (and references therein).

We have given special attention to the long-range dependence stable process $\xi$ and its properties because long-range dependence is currently a subject of much interest (see e.g. \cite{DOT}, \cite{H1}, \cite{H2}, \cite{S}, \cite{T} for discussions and literature), hence it is worthwhile to study different types of stochastic models where it appears. Other types of long-range dependence processes have been found recently (e.g. \cite{CS}, \cite{GNR}, \cite{HJ}, \cite{HV}, \cite{KT}, \cite{MY}, \cite{LT}, \cite{PTL}), in particular in models involving heavy-tailed distributions.

The following notation is used in the paper.

${\cal S}(\erre^d)$: space of $C^\infty$ rapidly decreasing functions on $\erre^d$.

${\cal S}'(\erre^d)$: space of tempered distributions (topological dual of 
${\cal S}(\erre^d)$).

$\langle\quad,\quad\rangle$: duality on ${\cal S}'(\erre^d)\times 
{\cal S}(\erre^d)$, 
or on ${\cal S}'(\erre^{d+1})\times {\cal S}(\erre^{d+1})$,
in particular, integral of a function with respect to a tempered measure.

$\Rightarrow_C$: weak convergence on the space of continuous functions 
$C([0,\tau], {\cal S}'(\erre^d))$  for each $\tau>0$. 


$\Rightarrow_f$: weak convergence of finite-dimensional distributions of 
${\cal S}'(\erre^d)$-valued processes.

$\Rightarrow_i$: integral convergence of ${\cal S}'(\erre^d)$-valued processes, i.e., $X_T\Rightarrow_i X$ if, for any $\tau>0$, the ${\cal S}'(\erre^{d+1})$-valued random variables $\widetilde{X}_T$ converge in law to 
$\widetilde{X}$ as $T\to\infty$, where $\widetilde{X}$ (and, analogously $\widetilde{X}_T$) is defined as a space-time random field by
\begin{equation}
\label{eq:1.11}
\langle\widetilde{X},\Phi\rangle=\int^\tau_0\langle X(t),\Phi(\cdot,t)\rangle dt,\quad\Phi\in{\cal S}(\erre^{d+1}).
\end{equation}

$\Rightarrow_{f,i}:\,\,\Rightarrow_f$ and $\Rightarrow_i$ together.

Recall
that in general $\Rightarrow_f$ and $\Rightarrow_i$ do not imply each other, and either one of them, together with tightness of $\{\langle X_T,\varphi\rangle\}_{T\ge 1}$ in $C([0,\tau],\erre)$, $\tau>0$, $\varphi\in{\cal S}(\erre^d)$, implies $\Rightarrow_C$ \cite{BGR}.

The transition probability density, the semigroup, and the potential operator of the standard symmetric $\alpha$-stable L\'evy process on $\erre^d$ are denoted respectively by $p_t(x), {\cal T}_t$ 
(i.e.,
${\cal T}_t\varphi=p_t*\varphi$)
and (for $d>\alpha$)
\begin{equation}
\label{eq:1.12}
G\varphi(x)=\int^\infty_0{\cal T}_t\varphi(x)dt=C_{\alpha,d}\int_{\errita}
\frac{\varphi(y)}{|x-y|^{d-\alpha}}dy,
\end{equation}
where
\begin{equation}
\label{eq:1.13}
C_{\alpha,d}=\frac{\Gamma(\frac{d-\alpha}{2})}{2^\alpha\pi^{d/2}\Gamma(\frac{\alpha}{2})}.
\end{equation}

Generic constants are written, $C, C_i$, with possible dependencies in parenthesis.

Section 2 contains the results, and Section 3 the proofs.

\vglue1cm
\noindent
{\bf 2. Results}
\setcounter{section}{2}
\setcounter{equation}{0}
\vglue.5cm

Given $\beta\in (0,1]$, let $M$ be an independently scattered $(1+\beta)$-stable measure on $\erre^{d+1}$ with control measure $\lambda_{d+1}$ (Lebesgue measure) and skewness intensity 1, i.e., for each Borel set $A\subset \erre^{d+1}$ 
such that 
$0< \lambda_{d+1}(A)<\infty, M(A)$ is a $(1+\beta)$-stable random variable with characteristic function
$${\rm exp}\left\{-\lambda_{d+1}(A)|z|^{1+\beta}
\left(1-i({\rm sgn}\, z)\tan\frac{\pi}{2}(1+\beta)\right)\right\},\quad 
z\in \erre,$$
the values of $M$ are independent on disjoint sets, and $M$ is $\sigma$-additive a.s. (see \cite{ST}, Definition 3.3.1).

For $\alpha\in (0,2], \gamma\geq 0$, we define the process
\begin{equation}
\label{eq:2.1}
\xi_t=\int_{\errita^{d+1}}\left(\UNO_{[0,t]}(r)
\left(\int_{\errita^d}p_r(x-y)|y|^{-\gamma}dy\right)^{1/(1+\beta)}
\int^t_rp_{u-r}(x)du
\right)M(drdx),\quad t\geq 0,
\end{equation}
which is well defined provided that
\begin{equation}
\label{eq:2.2}
\int_{\errita^d}\int^t_0\int_{\errita^d}p_r(x-y)|y|^{-\gamma}dy\left(\int^t_rp_{u-r}(x)du\right)^{1+\beta}drdx<\infty,
\end{equation}
(see \cite{ST}). For $\gamma=0,\ \xi$ is the same as the process $\xi$ defined by (2.1)
 in \cite{BGT6}. We also recall the following process defined by (2.2) in 
\cite{BGT6}, 
\begin{equation}
\label{eq:2.3}
\zeta_t=\int_{\errita^{d+1}}\left(\UNO_{[0,t]}(r)p_r^{1/(1+\beta)}(x)
\int^t_rp_{u-r}(x)du\right)M(drdx),\quad t\geq 0,
\end{equation}
which is well defined if $d<\alpha (2+\beta)/(1+\beta)$.

We consider the $(d,\alpha,\beta,\gamma)$-branching particle system described in the Introduction with $X_T$ defined by (\ref{eq:1.7}). 
Recall that the initial Poisson intensity measure is $H_T\mu_\gamma$.
We formulate the results for low, critical and high dimensions separately, since, as mentioned in the Introduction, the qualitative behaviors of the limit processes are different in
 each one of these cases. In the theorems below $K$ is a positive number depending on $d,\alpha,\beta,\gamma,V$, which may vary from case
to case
 and it is possible to compute it explicitly.                          

The results for the  low dimensions are contained in the following theorem.
\begin{theorem}
\label{T:2.1}
(a) Assume  $\gamma<d$ and
\begin{equation}
\label{eq:2.4}
d<\alpha\frac{2+\beta}{\beta}-\frac{\gamma\vee \alpha}{\beta}.
\end{equation}
Then the process $\xi$ given by (\ref{eq:2.1}) is well defined, and for
\begin{equation}
\label{eq:2.5}
 F^{1+\beta}_T=H_TT^{2+\beta-(d\beta+\gamma)/\alpha}
\end{equation}
with $H_T \geq 1$ and
\begin{equation}
\label{eq:2.6}
\lim_{T\to\infty}\frac{T^{1-(d-\gamma)\beta/\alpha}}{H^\beta_T}=0,
\end{equation}
we have $X_T\Rightarrow_CK\lambda\xi$.

\noindent
(b) Let $\gamma\geq d$,
\begin{equation}
\label{eq:2.7}
d<\alpha \frac{2+\beta}{1+\beta},
\end{equation}
and put
\begin{equation}
\label{eq:2.8}
k(T)=\left\{\begin{array}{lll}
\log T&{\it if}&\gamma=d,\nonumber\\
1&{\it if}&\gamma>d.
\end{array}\right.
\end{equation}
Then for
\begin{equation}
\label{eq:2.9}
F^{1+\beta}_T=H_TT^{2+\beta-(1+\beta)d/\alpha}k(T)
\end{equation}
with
\begin{equation}
\label{eq:2.10}
\lim_{T\to\infty}\frac{T}{H^\beta_Tk(T)^\beta}=0,
\end{equation}
we have $X_T\Rightarrow_{f,i}K\lambda\zeta$, where $\zeta$ is defined by (\ref{eq:2.3}).
\end{theorem}

\begin{remark}
\label{R:2.2}
{\rm (a) For $\gamma=0$, Theorem 2.1(a) is the same as Theorem 2.2(a) in 
\cite{BGT6}. For $\gamma\leq \alpha$ and $\gamma<d$, the bound on the dimension remains the same, equal to $\alpha(1+\beta)/\beta$ (see (\ref{eq:2.4})), and for $\alpha<\gamma<d$, it changes continuously, tending to the threshold 
(\ref{eq:2.7}) as $\gamma\nearrow d$.}
\end{remark}
\vglue-.5cm
\noindent
(b) For $d$ satisfying (\ref{eq:2.4}) and additionally $d>{\alpha}/{\beta}+\gamma$, condition (\ref{eq:2.6}) holds with $H_T=1$, so in this case  high density is not needed, and the limit of $X_T$  is the same as for the high-density model.

\noindent
(c) The case $\gamma> d$ is included for completeness only, since it is contained in Theorem 2.7 of \cite{BGT6}, where a general finite intensity measure was considered. The same remark applies also to the theorems for critical  and high dimensions (Theorems 2.5 and 2.6 below).

\noindent
(d) Note that the limit process is the same (up to constant) for the infinite intensity measure $H_Tdx/(1+|x|^d)\quad  (\gamma=d)$ as for finite measures.

\noindent
(e) In  Theorem 2.1(b) we consider convergence $\Rightarrow_{f,i}$ only. We are sure that  functional convergence holds (in fact, for the case $\gamma>d$ this was proved in \cite{BGT6}), but, as stated in the Introduction, we are mainly interested in the identification of limits and we do not attempt to give convergence results in the strongest forms. The same  applies to the theorems that follow. 

In the next proposition we gather some basic properties of the process $\xi$ defined by (\ref{eq:2.1}), in particular its
  long-range dependence property. (The process $\zeta$ was discussed in 
\cite{BGT6}). In \cite{BGT3} we introduced a way of measuring  long-range dependence in terms of the {\it dependence exponent}, defined by 
\begin{equation}
\label{eq:2.11}
\kappa=\inf_{z_1,z_2\in\errita}\inf_{0\leq u<v< s< t}
\sup\{\theta>0: D_T(z_1,z_2;u,v,s,t)=o(T^{-\theta})\,\,{\rm as}\,\, T\to\infty\},
\end{equation}
where
\begin{eqnarray}
\lefteqn{D_T(z_1,z_2;u,v,s,t)}\nonumber\\
\label{eq:2.12}
&=&|\log E e^{i(z_1(\xi_v-\xi_u)+z_2(\xi_{T+t}-\xi_{T+s})}-\log Ee^{iz_1(\xi_v-\xi_u)}-\log Ee^{iz_2(\xi_{T+t}-\xi_{T+s})}|.
\end{eqnarray}
(see also \cite{RZ}).
\begin{proposition}
\label{P:2.3}
Assume $\gamma<d$ and (\ref{eq:2.4}). Then

\noindent
(a) $\xi$ is $(1+\beta)$-stable, totally skewed to the right if $\beta<1$.

\noindent
(b) $\xi$ is self-similar with index 
$(2+\beta-(d\beta-\gamma)/\alpha)/(1+\beta)$.

\noindent
(c) $\xi$ has continuous paths.

\noindent
(d) $\xi$ has the long-range dependence property with dependence exponent

\begin{equation}
\label{eq:2.13}
\kappa=\left\{\begin{array}{l}
\frac{d}{\alpha}\quad\hbox{\it if either}\quad \alpha=2,\quad{\it or}\quad\alpha<2\quad{\it and}\quad \beta>\frac{d-\gamma}{d+\alpha},\nonumber\\
\frac{d}{\alpha}(1+\beta-\frac{d-\gamma}{d+\alpha})\quad{\it if}\quad\alpha<2\quad{\it and}\quad \beta\leq \frac{d-\gamma}{d+\alpha}.
\end{array}\right.
\end{equation}
\end{proposition}

\begin{remark}
\label{R:2.4}
{\rm (a) Here, as in the case $\gamma=0$ (Theorem 2.7 of \cite{BGT3}), the 
intriguing phenomenon of two long-range dependence regimes occurs 
for $\alpha<2$. It seems also interesting to note 
that putting formally $\gamma\geq d$ in 
(\ref{eq:2.13}) we obtain $\kappa={d}/{\alpha}$ (with no change of regime), which is indeed the dependence exponent of the process $\zeta$ (Proposition 2.9 of \cite{BGT6}). On the other hand, the process $\zeta$ itself is not obtained from $\xi$ by putting $\gamma\geq d$.

(b) If $\gamma=0$ and $\beta=1$, then $\xi$ is the sub-fractional Brownian motion (multiplied by a constant) considered in \cite{BGT}, \cite{BGT1}.}
\end{remark}

We now turn to the critical dimensions, i.e.,  the cases where the inequalities in (\ref{eq:2.4}) and (\ref{eq:2.7}) are replaced by equalities. It turns out that in spite of different conditions on the normings, the limits have always the same form as for finite intensity measure, with the only exception of the case given in  Theorem \ref{T:2.5}(a) below.

\begin{theorem}
\label{T:2.5}
(a) Assume $\gamma<d, \gamma<\alpha$,
\begin{equation}
\label{eq:2.14}
d=\alpha\frac{1+\beta}{\beta}
\end{equation}
and
\begin{equation}
\label{eq:2.15}
F^{1+\beta}_T=H_TT^{1-\gamma/\alpha}\log T,
\end{equation}
with $H_T\geq 1$. Then $X_T\Rightarrow_{f,i}K \lambda\eta$, where $\eta$ is a $(1+\beta)$-stable process with independent, non-stationary increments (for $\gamma>0$) whose laws are determined by
$$Ee^{iz(\eta_t-\eta_s)}={\rm exp}
\left\{-(t^{1-\gamma/\alpha}-s^{1-\gamma/\alpha})|z|^{1+\beta}
\left(1-i({\rm sgn} z)\tan\frac{\pi}{2}(1+\beta)
\right)\right\},\quad z\in\erre,\quad t\geq s\geq 0,$$%
$\eta_0=0$.

\noindent
(b) In all the remaining critical cases, i.e.,

\noindent
(i) $\gamma=\alpha, \gamma<d$ with $d$ satisfying (\ref{eq:2.14}),
\begin{equation}
\label{eq:2.16}
F^{1+\beta}_T=H_T(\log T)^2,
\end{equation}
and
\begin{equation}
\label{eq:2.17}
\lim_{T\to\infty}\frac{(\log T)^{1-\beta}}{H^\beta_T}=0,
\end{equation}
(ii) $\alpha<\gamma<d$ with
\begin{eqnarray}
\label{eq:2.18}
d=\alpha\frac{2+\beta}{\beta}-\frac{\gamma}{\beta},\\
\label{eq:2.19}
F^{1+\beta}_T=H_T\log T,
\end{eqnarray}
and
\begin{equation}
\label{eq:2.20}
\lim_{T\to\infty}\frac{T^{(1+\beta)(\gamma/\alpha-1)}}{H^\beta_T(\log T)^{1+\beta}}=0,
\end{equation}

\noindent
(iii)
\begin{equation}
\label{eq:2.21}
\gamma=d=\alpha\frac{2+\beta}{1+\beta},
\end{equation}
with $F_T$ satisfying (\ref{eq:2.16}) and
\begin{equation}
\label{eq:2.22}
\lim_{T\to\infty}\frac{T}{H^\beta_T(\log T)^{1+2\beta}}=0,
\end{equation}
(iv) $\gamma>d =\alpha ({2+\beta})/({1+\beta}), F^{1+\beta}_T=H_T \log T$ and $\lim_{T\to\infty}T H^{-\beta}_T=0$,

\noindent
we have $X_T\Rightarrow_{f,i}K\lambda\vartheta$, where $\vartheta$ is a real process such that $\vartheta_0=0$ and for $t>0, \vartheta_t=\vartheta_1=(1+\beta)$-stable random variable totally skewed to the right, i.e.,
$$Ee^{iz\vartheta_1}={\rm exp}\left\{-|z|^{1+\beta}\left(1-i({\rm sgn} z)\tan 
\frac{\pi}{2}(1+\beta)\right)\right\}.$$
\end{theorem}

It remains to consider the high dimensions.

\begin{theorem}
\label{T:2.6}
(a) Assume $\gamma<d, \gamma<\alpha$,
\begin{equation}
\label{eq:2.23}
d>\alpha\frac{1+\beta}{\beta},
\end{equation}
and
\begin{equation}
\label{eq:2.24}
F^{1+\beta}_T=H_TT^{1-{\gamma}/{\alpha}},
\end{equation}
with $H_T\geq 1$. Then $X_T\Rightarrow_{f,i}X$, where $X$ is an ${\cal S}'(\erre^d)$-valued $(1+\beta)$-stable process with independent, non-stationary increments (for $\gamma>0$) determined by
\begin{eqnarray}
\lefteqn{Ee^{i\langle X_t-X_s,\varphi\rangle}}\nonumber\\
&&={\rm exp}
 \Big\{-K(t^{1-\gamma/\alpha}-
s^{1-\gamma/\alpha})\int_{\errita^d}\Big(V|G\varphi(x)|^{1+\beta}(1-i({\rm sgn} G\varphi(x))\tan\frac{\pi}{2}(1+\beta))\nonumber\\
&&
+\,\,2c_\beta\varphi(x)G\varphi(x)\Big)dx\Big\},
\phantom{dddddddddddddddddddddddd}
\varphi\in{\cal S}(\erre^d),\quad t\geq 
s\geq 0,\label{eq:2.25}
\end{eqnarray}
$X_0=0$, where
\begin{equation}
\label{eq:2.26}
c_\beta=\left\{\begin{array}{lll}
0&{\it if}&0<\beta<1,\\
1&{\it if}&\beta=1,\nonumber
\end{array}\right.
\end{equation}
and $G$ is defined by (\ref{eq:1.12}).

\noindent
(b) Assume $\gamma<d,\gamma=\alpha$, and $d$ satisfying (\ref{eq:2.23}) with
\begin{equation}
\label{eq:2.27}
F^{1+\beta}_T=H_T\log T,
\end{equation}
and $H_T\geq 1$. Then $X_T\Rightarrow_{f,i}X$, where $X$ is an ${\cal S}'(\erre^d)$-valued process such that $X_0=0$, and for $t>0, X_t=X_1=(1+\beta)$-stable random variable determined by
\begin{eqnarray}
Ee^{i\langle X_1,\varphi\rangle}
&=&{\rm exp}\biggl\{-K\int_{\errita^d}\biggr.
\left(V|G\varphi(x)|^{1+\beta}(1-i({\rm sgn}
 G\varphi(x))\tan\frac{\pi}{2}(1+\beta)\right.)\nonumber\\
\label{eq:2.28}
&+&2c_\beta\varphi(x)G\varphi(x)\Big)dx\biggr\},\quad \biggl.
\phantom{dddddddddddddddddd}
\varphi\in{\cal S}(\erre^d),
\end{eqnarray}
with $c_\beta$ given by (\ref{eq:2.26}).

\noindent
(c) Assume $\gamma>\alpha$,
\begin{eqnarray}
\label{eq:2.29}
d>\alpha\frac{2+\beta}{\beta}&-&\frac{\gamma\wedge d}{\beta},\\
\label{eq:2.30}
F^{1+\beta}_T&=&H_T,
\end{eqnarray}
and
\begin{equation}
\label{eq:2.31}
\lim_{T\to\infty} TH^{-\beta}_T=0.
\end{equation}
Then $X_T\Rightarrow_{f,i}X$, where $X$ is an ${\cal S}'(\erre^d)$-valued process such that $X_0=0$, and for $t>0, X_t=X_1=(1+\beta)$-stable random variable determined by
\begin{eqnarray}
\lefteqn{
\kern-1cm Ee^{i\langle X_1,\varphi\rangle}={\rm exp}
\left\{-K\int_{\errita^d}\right.\left(
V|G\varphi(x)|^{1+\beta}(1-i({\rm sgn} G\varphi(x))\tan\frac{\pi}{2}(1+\beta)
\right.)\nonumber}\\
\label{eq:2.32}
&&+\,\,2c_\beta\varphi(x)G\biggl.\biggl.\varphi(x)\biggr)G\mu_\gamma(dx)\biggr\},\quad\phantom{ddddddddddd} \varphi\in 
{\cal S}(\erre^d),
\end{eqnarray}
with $c_\beta$ given by (\ref{eq:2.26}).
\end{theorem}

\begin{remark}
\label{R:2.7}
{\rm (a) As in all the cases studied previously 
\cite{BGT3}, \cite{BGT4}, \cite{BGT6}, we observe the same phenomenon that in low dimensions the limit processes are continuous with a simple spatial structure and a complicated temporal structure (with long-range dependence), while in high dimensions they are truly ${\cal S}'(\erre^d)$-valued with independent increments, and not necessarily continuous.}
\end{remark}
\vglue-.4cm
\noindent
(b) For low dimensions the forms of the limits depend on the relation between $d$ and $\gamma$ only, whereas for critical and high dimensions only the relationship between $\alpha$ and $\gamma$ is relevant. More precisely, in critical dimensions we have different forms of the limits for $\gamma<\alpha$ and $\gamma\geq \alpha$, and in high dimensions the forms are different for $\gamma<\alpha, \gamma=\alpha$ and $ \gamma>\alpha$. In the  case $\gamma>\alpha$ even the 
normings are the same, depending only on $\beta$.

\noindent
(c) For $\beta=1$ the limits are centered Gaussian.
In high dimensions there is no continuous transition between the cases $\beta<1$ and $\beta=1$; in the latter case an additional term appears. 
The coefficient $c_\beta$ defined in (\ref{eq:2.26}) was introduced in order to present the results in unified forms.

\noindent
(d) We have assumed that the initial intensity  measure is determined 
by $\mu_\gamma$ of the form
(\ref{eq:1.8}). It will be  clear that the same results are obtained for the measure 
$\mu_\gamma(dx)=|x|^{-\gamma}dx$ if $d>\gamma$. Analogously as in the non-branching case \cite{BGT5}, other generalizations are also possible.

Let us look further into the need for high density  (i.e., to assume 
$H_T\to\infty$) and the  question of a.s. local extinction.  
For $\gamma>d$  the Poisson intensity measure is finite,  
so there is only a finite number of particles at time $t=0$ and the system becomes globally extinct in 
finite time a.s. due to the criticality of the branching.
%
%
%
 Also, for $\gamma\wedge d>\alpha$ it is not difficult to see that the total occupation time $\int^\infty_0N_sds$ is finite a.s. on bounded sets (see 
\cite{BGT5}, Proposition~2.1, because $EN_s$ is  the same for the systems with and without branching), so high density is also necessary.  We have a more delicate situation in the remaining cases where the threshold is $d=\alpha/\beta+\gamma$. Concerning  extinction,
the situation is completely clear for $\alpha=2$. In Theorem 2.8 below we state
 that for $\alpha=2$ and $d< 2 /\beta+\gamma$ there is a.s. local extinction,
 hence the total occupation time of any bounded set is finite a.s. We conjecture that the same is true for $d<\alpha /\beta +\gamma$ if $\alpha<2$, but we have not been able to prove it.

\begin{theorem}
\label{T:2.8}
Assume $\alpha=2$. If $d<2/\beta+\gamma$, then for each bounded Borel set $A$,
$$P[\hbox{\it there exists}\quad \tau_A<\infty\quad\hbox{\it such that}\quad N_t(A)=0
\quad\hbox{\it for all}\quad t\geq\tau_A]=1.$$
\end{theorem}

The proof of this theorem relies on Iscoe's  a.s.\
local extinction result for the 
superprocess   \cite{I2}, 
 by showing that in general (i.e., for $0<\alpha\le 2$) a.s.\ local finiteness of the total occupation time of the $(d,\alpha,\beta,\gamma)$-superprocess implies a.s.\ local extinction of the $(d,\alpha,\beta,\gamma)$-branching particle system. On the other hand, as explained in the Introduction, for $\alpha<2$ the a.s.\ local extinction for the superprocess cannot occur, and we do not know how to prove directly the a.s.\ local finiteness of its total occupation time.

The next 
ergodic-type result, which is a direct generalization of 
\cite{Ta}, shows that $\alpha/\beta+\gamma$ is indeed a natural threshold.

\begin{proposition}
\label{P:2.9}
Assume
\begin{equation}
\label{eq:2.33}
\gamma<\alpha,\quad d=\frac{\alpha}{\beta}+\gamma,\quad
 F_T=T^{1-{\gamma}/{\alpha}},
\end{equation}
and denote
\begin{equation}
\label{eq:2.34}
Z_T(t)=\frac{1}{F_T}\int^{Tt}_0N_sds,\quad t\geq 0.
\end{equation}
Then $Z_T\Rightarrow_C\lambda\xi$, 
where $\xi$ is a real non-negative process with Laplace transform
\begin{equation}
\label{eq:2.35}
E{\rm exp}\{-\theta_1\xi_{t_1}-\cdots -\theta_n\xi_{t_n}\}={\rm exp}
\left\{-\int_{\errita^d}v(x,\tau)|x|^{-\gamma}dx\right\},
\end{equation}
for 
any $\tau>0$, where $\theta_1,\ldots,\theta_n\geq 0, 0\leq t_1<\cdots <t_n\leq \tau$, and 
$v(x,t)$ is the unique non-negative solution of the equation
\begin{equation}
\label{eq:2.36}
v(x,t)=\int^t_0p_{t-s}(x)\psi(\tau-s)ds-\frac{V}{1+\beta}\int^t_0{\cal T}_{t-s}v^{1+\beta}(\cdot,s)(x)ds,
\end{equation}
with
\begin{equation}
\label{eq:2.37}
\psi(s)=\sum^n_{k=1}\theta_k\UNO_{[0,t_k]}(s).
\end{equation}
\end{proposition}

To complete the discussion on a.s.\ local extinction we formulate a
 corollary which follows immediately from our results but which, nevertheless, seems worth stating explicitly.

\begin{corollary}
\label{C:2.10}
If $\gamma<\alpha$ and $d\ge \alpha/\beta+\gamma$, then the $(d,\alpha,\beta,\gamma)$-branching particle system does not have the a.s.\ local extinction property.
\end{corollary}

Indeed, for $d>\alpha/\beta+\gamma$, from Theorems 2.1, 2.5 and 2.6 it follows that one may take $H_T\equiv 1$. By Proposition 2.1 of \cite{BGT5}, $E\int_0^T\langle N_s,\varphi\rangle ds$ (for $\int_{\errita^d}\varphi(x)dx\ne 0$) is of 
larger order than $F_T$ as $T\to\infty$ hence, by (\ref{eq:1.7}), for any  
bounded Borel set $ A$, $\int_0^\infty N_s(A)ds=\infty$ a.s., which excludes 
a.s. local extinction. For $d=\alpha/\beta+\gamma$ the result follows immediately from Proposition 2.9.

We end with the results for the superprocess.
\vglue.5cm
\noindent
{\bf Theorem 2.11} {\it Let $Y$ be the $(d,\alpha,\beta,\gamma)$-superprocess and $X_T$  its
 occupation time fluctuation process defined by (1.10). Then the limit results for $X_T$ as $T\to\infty$ are the same as those in Theorems 2.1, 2.5 and 2.6, with the same normings, and $c_\beta=0$ in all cases in Theorem 2.6.}

\vglue1cm
\noindent
{\bf 3. Proofs}
\vglue.5cm
\noindent
{\bf 3.1 Scheme of proofs}
\setcounter{section}{3}
\setcounter{equation}{0}
\vglue.5cm
The proofs of Theorems 2.1, 2.5 and 2.6 follow the general scheme presented in \cite{BGT6}. For completeness we recall the main steps. 

As explained in \cite{BGT3,BGT4,BGT6}, in order to prove convergence $\Rightarrow_i$ it suffices to show
\begin{equation}
\label{eq:3.1}
\lim_{T\to\infty} Ee^{-\langle\widetilde{X}_T,\Phi\rangle}=Ee^{-\langle\widetilde{X},\Phi\rangle}
\end{equation}
for each $\Phi\in {\cal S}(\erre^{d+1}),\Phi\geq 0$, where $X$ is the corresponding limit process and $\widetilde{X}_T,\widetilde{X}$ are defined by 
(\ref{eq:1.11}).
To prove convergence $\Rightarrow_C$ according to the space-time approach \cite{BGR} it is enough to show additionally that the family $\{\langle X_T,\varphi\rangle\}_{T\geq 1}$ is tight in
 $C((0,\tau],\erre),\varphi\in {\cal S}(\erre^d), \tau>0$. Without loss of generality we may fix $\tau=1$ (see (\ref{eq:1.11})).
To simplify slightly the calculations we consider $\Phi$ of the form
$$\Phi(x,t)=\varphi\otimes\psi(x,t)=\varphi(x)\psi(t),\,\,\varphi\in
{\cal S}(\erre^d), \psi\in{\cal S}(\erre),\varphi,\psi\geq 0.$$
Denote 
\begin{equation}
\label{eq:3.2}
\varphi_T=\frac{1}{F_T}\varphi,\,\, \chi(t)=\int^1_t\psi(s)ds,\,\, 
\chi_T(t)=\chi\left(\frac{t}{T}\right).
\end{equation}

We define
\begin{equation}
\label{eq:3.3}
v_T(x,t)=1-E{\rm exp}\left\{-\int^t_0\langle N^x_r,\varphi_T\rangle\chi_T(T-t+r)dr\right\},\,\,0\leq t\leq T,
\end{equation}
where $N^x$ is the empirical process of the branching system started from a single particle at $x$. The following equation for $v_T$ was derived in \cite{BGT3} (formula (\ref{eq:3.8}), see also \cite{BGT1}) by means of the Feynman-Kac formula:
\begin{equation}
\label{eq:3.4}
v_T(x,t)=\int^t_0{\cal T}_{t-u}\left[\varphi_T\chi_T(T-u)(1-v_T(\cdot,u))
-\frac{V}{1+\beta}v^{1+\beta}_T(\cdot,u)\right](x)du,
\,\,\,0\leq t\leq T.
\end{equation}
This equation (with $T=1$) is the space-time version of equation (1.4). It is the log-Laplace equation for $\widetilde{L}$ (as in (1.11)), where $L$ is the occupation time $L_t=\int^t_0N_sds$. Formulas
(\ref{eq:3.3}) and (\ref{eq:3.4}) imply
\begin{equation}
\label{eq:3.5}
0\leq v_T\leq 1,\quad v_T(x,t)\leq \int^t_0{\cal T}_{t-u}\varphi_T(x)\chi_T(T-u)du.
\end{equation}

For brevity we denote
\begin{equation}
\label{eq:3.6}
\nu_T(dx)=H_T\mu_\gamma(dx)=\frac{H_T}{1+|x|^\gamma}dx.
\end{equation}
By the Poisson property and (\ref{eq:3.4}) we have
\begin{eqnarray}
\label{eq:3.7}
Ee^{-\langle\widetilde{X}_T,\varphi\otimes\psi\rangle}&=&{\rm exp}
\left\{-\int_{\errita^d}v_T(x,T)\nu_T(dx)+\int_{\errita^d}\int^T_0{\cal T}_u\varphi_T(x)\chi_T(u)du\nu_T(dx)\right\}\\
\label{eq:3.8}
&=&{\rm exp}\left\{\frac{V}{1+\beta}I_1(T)+I_2(T)-\frac{V}{1+\beta}I_3(T)
\right\},
\end{eqnarray}
where
\begin{eqnarray}
\label{eq:3.9}
I_1(T)&=&\int_{\errita^d}\int^T_0{\cal T}_{T-s}
\left[\left(\int^s_0{\cal T}_{s-u}\varphi_T\chi_T(T-u)du\right)^{1+\beta}\right](x)ds\nu_T(dx),\\
\label{eq:3.10}
I_2(T)&=&\int_{\errita^d}\int^T_0{\cal T}_{T-s}(\varphi_T\chi_T(T-s)v_T(\cdot,s))(x)ds\nu_T(dx),\\
\label{eq:3.11}
I_3(T)&=&\int_{\errita^d}\int^T_0{\cal T}_{T-s}\left[\left(
\int^s_0{\cal T}_{s-u}\varphi_T\chi_T(T-u)du\right)^{1+\beta}-v^{1+\beta}_T(\cdot,s)\right](x)ds\nu_T(dx).
\end{eqnarray}
In the proofs of Theorems 2.1, 2.5 and in Theorem 2.6 for $\beta<1$ we show
\begin{equation}
\label{eq:3.12}
\lim_{T\to\infty}
{\rm exp}\left\{\frac{V}{1+\beta}I_1(T)\right\}=Ee^{-\langle\widetilde{X},
\varphi\otimes\psi\rangle},
\end{equation}
and
\begin{equation}
\label{eq:3.13}
\lim_{T\to\infty}I_2(T)=0,
\end{equation}
where (\ref{eq:3.13}) is obtained from
\begin{equation}
\label{eq:3.14}
I_2(T)\leq \frac{C}{F^2_T}\int_{\errita^d}\int^T_0\int^T_0{\cal T}_s(\varphi{\cal T}_u\varphi)(x)duds\nu_T(dx)
\end{equation}
(see (\ref{eq:3.5})).
In Theorem 2.6 for $\beta=1$ the limit of $I_2(T)$ is non-trivial and corresponds to the expressions involving $c_\beta$ (see (\ref{eq:2.25}), (\ref{eq:2.26}),
(\ref{eq:2.28}),(\ref{eq:2.32})). In all the cases
\begin{equation}
\label{eq:3.15}
\lim_{T\to\infty}I_3(T)=0
\end{equation}
By the argument in \cite{BGT6}, in order to prove (\ref{eq:3.15}) we show
\begin{equation}
\label{eq:3.16}
\lim_{T\to\infty}J_1(T)=0
\end{equation}
and
\begin{equation}
\label{eq:3.17}
\lim_{T\to\infty} J_2(T)=0,
\end{equation}
where
\begin{eqnarray}
J_1(T)&=&\int_{\errita^d}\int^T_0{\cal T}_{T-s}
\left[\left(\int^s_0{\cal T}_{s-u}\left(\varphi_T\int^u_0{\cal T}_{u-r}
\varphi_Tdr\right)du\right)^{1+\beta}\right](x)ds\nu_T(dx)\nonumber\\
\label{eq:3.18}
&\leq&\frac{1}{F^{2+2\beta}_T}\int_{\errita^d}\int^T_0{\cal T}_s
\left[\left(\int^T_0{\cal T}_u\left(\varphi\int^T_0{\cal T}_r\varphi 
dr\right)du\right)^{1+\beta}\right](x)ds\nu_T(dx),\\
J_2(T)&=&\int_{\errita^d}\int^T_0{\cal T}_{T-s}
\left[\left(\int^s_0{\cal T}_{s-u}
\left(\int^u_0{\cal T}_{u-r}\varphi_Tdr\right)^{1+\beta}du\right)^{1+\beta}
\right](x)ds\nu_T(dx)\nonumber\\
\label{eq:3.19}
&\leq&\frac{1}{F^{(1+\beta)(1+\beta)}_T}\int_{\errita^d}\int^T_0{\cal T}_s
\left[\left(\int^T_0{\cal T}_u\left(\int^T_0{\cal T}_r\varphi dr
\right)^{1+\beta}du\right)^{1+\beta}\right](x)ds\nu_T(dx).
\end{eqnarray}
We remark that the proof of (\ref{eq:3.15}) is the only place where the high density (with specific conditions on $H_T$) is required in some cases.

Finally, the $\Rightarrow_f$ convergence is obtained by an analogous argument as explained in the proof of Theorem 2.1 in \cite{BGT4}.

\vglue.5cm
\noindent
{\bf 3.2 Auxiliary estimates}
\vglue.5cm
Recall that the transition density $p_t$ of the standard $\alpha$-stable process has the self-similarity property
\begin{equation}
\label{eq:3.20}
p_{at}(x)=a^{-d/\alpha}p_t(a^{-1/\alpha}x),\quad x\in \erre^d,\,\, a>0,
\end{equation}
and it satisfies
\begin{equation}
\label{eq:3.21}
\frac{C_1}{1+|x|^{d+\alpha}}\leq p_1(x)\leq \frac{C_2}{1+|x|^{d+\alpha}},
\end{equation}
where the lower bound holds for $\alpha<2$.

Denote
\begin{equation}
\label{eq:3.22}
f(x)=\int^1_0p_s(x)ds.
\end{equation}
The following estimate can be easily deduced from (\ref{eq:3.20}) and 
(\ref{eq:3.21}):
\begin{eqnarray}
\label{eq:3.23}
f(x)&\leq& \frac{C}{|x|^{d+\alpha}},\\
f(x)&\leq& \left\{
\begin{array}{l}
C\quad{\rm if}\quad d<\alpha,\\
\label{eq:3.24}
C(1\vee \log |x|^{-1})\quad{\rm if}\quad d=\alpha,\\
{C}/|x|^{d-\alpha}\quad{\rm if}\quad d>\alpha. 
\end{array}\right.
\end{eqnarray}

We will also use the following elementary estimates: Let $0<a,b<d$. If $a+b>d,$ then
\begin{equation}
\label{eq:3.25}
\int_{\errita^d}\frac{1}{|x-y|^a|x|^b}dx\leq \frac{C}{|y|^{a+b-d}}.
\end{equation}
If $a+b=d$, then
\begin{equation}
\label{eq:3.26}
\int_{|x|\leq 1}\frac{1}{|x-y|^a|x|^b}dx\leq C(1\vee \log |y|^{-1}).
\end{equation}
If $a+b<d$, then
\begin{equation}
\label{eq:3.27}
\int_{|x|\leq 1}\frac{1}{|x-y|^a|x|^b}dx\leq C.
\end{equation}
Now, let $a>d, 0<b<d$, then
\begin{equation}
\label{eq:3.28}
\int_{\errita^d}\frac{1}{1+|x-y|^a}\frac{1}{|x|^b}dx\leq \frac{C}{|y|^b}.
\end{equation}
For $d>\gamma$, denote
\begin{equation}
\label{eq:3.29}
f_\gamma(y)=\int_{\errita^d}f(y-x)|x|^{-\gamma}dx,
\end{equation}
where $f$ is defined in (\ref{eq:3.22}). From the estimates above we obtain
\begin{equation}
\label{eq:3.30}
\sup_{|y|>1}|y|^\gamma f_\gamma(y)< \infty,
\end{equation}
and
\begin{equation}
\label{eq:3.31}
f_\gamma(y)\leq\left\{\begin{array}{l}
C\quad{\rm if}\quad \gamma<\alpha,\nonumber\\
C(1\vee\log |y|^{-1})\quad{\rm if}\quad \gamma=\alpha,\\
{C}/|y|^{\gamma-\alpha}\quad{\rm if}\quad \gamma>\alpha.\nonumber
\end{array}\right.
\end{equation}
\vglue.5cm
\noindent
{\bf 3.3 Proof of Theorem 2.1(a)}
\vglue.5cm
According to the scheme sketched above, in order to prove $\Rightarrow_{f,i}$ convergence we show (\ref{eq:3.1}). By (\ref{eq:3.7})-(\ref{eq:3.11})and 
(\ref{eq:2.1}) it is enough to prove (\ref{eq:3.12}), which amounts to 
\begin{equation}
\label{eq:3.32}
\lim_{T\to\infty}I_1(T)
=\int_{\errita^d}\int^1_0\int_{\errita^d}p_s(x-y)\left(\int^1_sp_{u-s}(y)
\chi(u)du\right)^{1+\beta}dyds\frac{dx}{|x|^\gamma}
\left(\int_{\errita^d}\varphi(z)dz\right)^{1+\beta},
\end{equation}
(see (\ref{eq:3.9})) and, additionally, (\ref{eq:3.13}) and (\ref{eq:3.15}). To simplify the notation we will carry out the proof for $\mu_\gamma$ of the form $\mu_\gamma(dx)=|x|^{-\gamma}dx$ instead of (\ref{eq:1.8}). It will be clear that in the present case $(d<\gamma)$ this will not 
affect the result.

By (\ref{eq:3.9}), (\ref{eq:3.2}), (\ref{eq:3.6}),
 the definition of ${\cal T}_t$, substituting $s'=1-s/T, u'=1-u/T$, we have
\begin{eqnarray}
\lefteqn{I_1(T)}\nonumber\\
\label{eq:3.33}
&=&\frac{T^{2+\beta}{H_T}}{F^{1+\beta}_T}\int_{\errita^d}\int^1_0\int_{\errita^d}p_{Ts}(x-y)\left(\int^1_s\int_{\errita^d}p_{T(u-s)}(y-z)\varphi(z)\chi(u)dzdu\right)^{1+\beta}|x|^{-\gamma}dydsdx.\nonumber\\
\end{eqnarray}

Denote
\begin{equation}
\label{eq:3.34}
\widetilde{\varphi}_T(x)=T^{d/\alpha}\varphi(T^{1/\alpha}x)
\end{equation}
and
\begin{equation}
\label{eq:3.35}
g_s(x)=\int^1_sp_{u-s}(x)\chi(u)du,\quad 0\leq s\leq 1.
\end{equation}
Observe that
\begin{equation}
\label{eq:3.36}
g_s\leq Cf,
\end{equation}
where $f$ is defined by (\ref{eq:3.22}). By (\ref{eq:3.20}) and (\ref{eq:2.5}),
making obvious spatial substitutions in (\ref{eq:3.33}),
 we obtain
\begin{equation}
\label{eq:3.37}
I_1(T)=\int_{\errita^d}\int^1_0\int_{\errita^d}p_s(x-y)(g_s*\widetilde{\varphi}_T(y))^{1+\beta}|x|^{-\gamma}dxdsdy.
\end{equation}
Note that if we consider the measure $\mu_\gamma$ of the form (\ref{eq:1.8}), then in (\ref{eq:3.37}) instead of $|x|^{-\gamma}$ we have $(T^{-\gamma/\alpha}+|x|^\gamma)^{-1}$. Since $g_s\in L^1(\erre^d)$, by (\ref{eq:3.34}) it is clear that $g_s * \widetilde{\varphi}_T(y)$ converges to $g_s(y)\int_{\errita^d}\varphi(z)dz$ almost everywhere in $y$. Hence, to prove (\ref{eq:3.32}) it remains to justify the passage to the limit under the integrals in (\ref{eq:3.37}). Denote
\begin{equation}
\label{eq:3.38}
h_T(y)=\int^1_0\int_{\errita^d}p_s(x-y)(g_s *\widetilde{\varphi}_T(y))^{1+\beta}|x|^{-\gamma}dxds.
\end{equation}
First we prove pointwise convergence of $h_T$, which amounts to showing that the integrand is majorized by an integrable function independent of $T$. Fix $y\neq 0$. We use (\ref{eq:3.36}) and observe that
\begin{eqnarray}
f*\widetilde{\varphi}_T(y)&=&\int_{|z|\leq {|y|}/{2}}f(y-z)\widetilde{\varphi}_T(z)dz+\int_{|z|>{|y|}/{2}}f(y-x)\frac{(T^{1/\alpha}|x|)^d\varphi(T^{d/\alpha}x)}{|x|^d}dx\nonumber\\
&\leq&f\left(\frac{y}{2}\right)\int_{|z|\leq{|y|}/{2}}\widetilde{\varphi}_T(z)dz+\frac{C}{|y|^d}\int_{|z|>{|y|}/{2}}f(y-x)dz\nonumber\\
\label{eq:3.39}
&\leq&f\left(\frac{y}{2}\right)\int_{\errita^d}\varphi(z)dz+\frac{C}{|y|^d},
\end{eqnarray}
by the unimodal property of the $\alpha$-stable density and since $\varphi\in{\cal S}(\erre^d)$.
We conclude by noting that
$$\int_{\errita^d}\int^1_0p_s(x-y)|x|^{-\gamma}dsdx<\infty\quad{\rm for}\quad y\neq 0,
$$
by (\ref{eq:3.31}).

Since (see (\ref{eq:3.29}))
\begin{equation}
\label{eq:3.40}
h_T(y)\leq \, f_\gamma(y)((f*\widetilde{\varphi}_T)(y))^{1+\beta},
\end{equation}
to prove convergence of $I_1(T)$ it suffices to show that the right-hand side of (\ref{eq:3.40})
 (denoted by $h^*_T)$ converges in $L^1(\erre^d)$ as $T\to\infty$.

If $\gamma<\alpha$, then $f_\gamma$ is bounded by (\ref{eq:3.31}), and the assumption (\ref{eq:2.4}) implies that
\begin{equation}
\label{eq:3.41}
f\in L^{1+\beta}(\erre^d),
\end{equation}
so $(f*\widetilde{\varphi}_T)^{1+\beta}$ converges in $L^1(\erre^d)$.

Next assume $\gamma\geq \alpha$. It is easily seen that $h^*_T\UNO_{\{|y|\geq 1\}}$ converges in $L^1(\erre^d)$, by (\ref{eq:3.41}) and (\ref{eq:3.30}). To prove that $h^*_T(y)\UNO_{\{|y|<1\}}$ converges in $L^1(\erre^d)$ too, it suffices to find $p,q>1, {1}/{p}+{1}/{q}=1$, such that
\begin{equation}
\label{eq:3.42}
f_\gamma(y)\UNO_{\{|y|<1\}}\in L^p(\erre^d)
\end{equation}
and
\begin{equation}
\label{eq:3.43}
f^{1+\beta}\in L^q(\erre^d).
\end{equation}

If $\gamma=\alpha$, then (\ref{eq:3.31}) implies that (\ref{eq:3.42})
 holds for any $p>1$, and by (\ref{eq:2.4}) it is clear that (\ref{eq:3.43}) is satisfied for $q$ sufficiently close to 1.

If $\gamma>\alpha$, condition (\ref{eq:2.4}) is equivalent to 
$$\frac{\gamma-\alpha}{d}+\frac{(1+\beta)(d-\alpha)}{d}<1,$$
so we can take $p$ and $q$ such that ${1}/{p}+{1}/{q}=1$,
$$\frac{1}{p}>\frac{\gamma-\alpha}{d}\quad{\rm and}\quad\frac{1}{q}>\frac{(1+\beta)(d-\alpha)}{d}.$$
For such $p$ and $q$ we have (\ref{eq:3.42}) and (\ref{eq:3.43}) by 
(\ref{eq:3.31}), (\ref{eq:3.23}) and (\ref{eq:3.24}).

This completes the proof of (\ref{eq:3.32}).

We proceed to the proof of (\ref{eq:3.13}). By (\ref{eq:3.14}), applying the same substitutions as for $I_1(T)$ and using the notation (\ref{eq:3.34}) we have
\begin{equation}
\label{eq:3.44}
I_2(T)\leq C\frac{H_TT^{2-d/\alpha-\gamma/\alpha}}{F^2_T}\int_{\errita^d}f_\gamma(y)\widetilde{\varphi}_T(y)(f *\widetilde{\varphi}_T)(y)dy.
\end{equation}
Assume $\gamma<\alpha$. By (\ref{eq:3.31}),
\begin{eqnarray*}
I_2(T)&\leq&C_1\frac{H_TT^{2-d/\alpha-\gamma/\alpha}}{F^2_T}||\widetilde{\varphi}_T(f*\widetilde{\varphi}_T)||_1\\
&\leq&C_1\frac{H_TT^{2-d/\alpha-\gamma/\alpha}}{F^2_T}||\varphi||_1||f||_{1+\beta}||\widetilde{\varphi}_T||_{\frac{1+\beta}{\beta}}\to 0
\end{eqnarray*}
by (\ref{eq:2.4}), and since
\begin{equation}
\label{eq:3.45}
||\widetilde{\varphi}_T||_p=||\varphi||_pT^{(d/\alpha)(p-1)/p}\quad{\rm for}\quad p\geq 1.
\end{equation}

Next, let $\gamma\geq \alpha$. By (\ref{eq:3.44})
 and (\ref{eq:3.30}) we have 
\begin{equation}
\label{eq:3.46}
I_2(T)\leq C_2(I'_2(T)+I^{''}_2(T)),
\end{equation}
where
\begin{equation}
\label{eq:3.47}
I'_2(T)=\frac{H_TT^{2-d/\alpha-\gamma/\alpha}}{F^2_T}||\widetilde{\varphi}_T(f*\widetilde{\varphi}_T)||_1,
\end{equation}
and
\begin{equation}
\label{eq:3.48}
I^{''}_2(T)=\frac{H_TT^{2-d/\alpha-\gamma/\alpha}}{F^2_T}\int_{|y|\leq 1}f_\gamma(y)\widetilde{\varphi}_T(y)(f*\widetilde{\varphi}_T)(y)dy.
\end{equation}
By the H\"older inequality,
\begin{equation}
\label{eq:3.49}
I'_2(T)\leq \frac{H_TT^{2-d/\alpha-\gamma/\alpha}}{F^2_T}||\widetilde{\varphi}_T||_p||f||_q||\varphi||_1
\end{equation}
for any $p,q\geq 1, 1/p+1/q=1$. If $1/q>(d-\alpha)/d$, then 
$||f||_q<\infty$, and if $1/q$ is sufficiently close to $(d-\alpha)/d$, then by (\ref{eq:3.45}), 
(\ref{eq:2.5}) and (\ref{eq:2.4}) the right-hand side of (\ref{eq:3.49}) converges to $0$ as $T\to\infty$.

We estimate $I^{''}_2(T)$ using the generalized H\"older inequality
$$I^{''}_2(T)\leq \frac{H_TT^{2-d/\alpha-\gamma/\alpha}}{F^2_T}||f_\gamma\UNO_{\{|\cdot|\leq 1\}}||_r||\widetilde{\varphi}_T||_p||f||_q||\varphi||_1$$
for $r,p,q\geq 1, 1/p+1/r+1/q=1$. We take $r,q$ such that $1/r>(\gamma-\alpha)/d$ (then the $r$-norm will be finite by (\ref{eq:3.31})) and $1/q>(d-\alpha)/d$. By (\ref{eq:2.4}), it is easily seen that one can choose $r,q$ as above and $p=\alpha/(2\alpha-\gamma)+\varepsilon$ for $\varepsilon>0$ sufficiently small (note that by (\ref{eq:2.4}), $2\alpha>\gamma$).
Then by (\ref{eq:3.45}) and (\ref{eq:2.5}) we obtain that $I^{''}_2(T)\to 0$ as $T\to\infty$. Thus, we have proved (\ref{eq:3.13}).

According to the general scheme, in order to obtain (\ref{eq:3.15}) it suffices to show (\ref{eq:3.16}) and (\ref{eq:3.17}). The proofs are quite similar to the argument presented above, therefore we omit the proof of (\ref{eq:3.16}) and we give an outline of the proof of (\ref{eq:3.17}), since this is the only place where the condition (\ref{eq:2.6}) is needed.

By (\ref{eq:3.19}), (\ref{eq:2.5}) and the usual substitutions we have
\begin{equation}
\label{eq:3.50}
J_2(T)\leq C\frac{T^{1-(d/\alpha)\beta+(\gamma/\alpha)\beta}}{H^\beta_T}R(T),
\end{equation}
where
\begin{equation}
\label{eq:3.51}
R(T)=\int_{\errita^d}f_\gamma(y)(f*(f*\widetilde{\varphi}_T)^{1+\beta})^{1+\beta}(y)dy.
\end{equation}
By (\ref{eq:2.6}), to prove (\ref{eq:3.17}) it remains to show that
\begin{equation}
\label{eq:3.52}
\sup_TR(T)<\infty.
\end{equation}

If $\gamma<\alpha$ then, by (\ref{eq:3.31}),
\begin{equation}
\label{eq:3.53}
R(T)\leq C_1||f*(f*\widetilde{\varphi}_T)^{1+\beta}||^{1+\beta}_{1+\beta}\leq C_1||f||^{1+\beta}_{1+\beta}||f||^{(1+\beta)(1+\beta)}_{1+\beta}||\varphi||^{(1+\beta)(1+\beta)}_1<\infty,
\end{equation}
by (\ref{eq:3.41}).

If $\gamma\geq \alpha$, then we write
$$R(T)=\int_{|y|>1}\ldots +\int_{|y|\leq 1}\ldots.$$
By (\ref{eq:3.30}) the first integral can be estimated as in (\ref{eq:3.53}), and the second one is bounded by
$$||f_\gamma\UNO_{\{|\cdot|\leq 1\}}||_p||f||^{1+\beta}_{q(1+\beta)}||f||^{(1+\beta)(1+\beta)}_{1+\beta}||\varphi||^{(1+\beta)(1+\beta)}_1,$$
where $1/p+1/q=1$. We already know that there exist such $p$ and $q$ that this expression is finite (see (\ref{eq:3.42}) and (\ref{eq:3.43})).

We have thus established the convergence
$$X_T\Rightarrow_{f,i}K\lambda\xi.$$

In order to obtain $\Rightarrow_C$ convergence it suffices to show that the family $\{\langle X_T,\varphi\rangle\}_{T\geq 1}$ is tight in $C([0,1],\erre)$ for any $\varphi\in {\cal S}(\erre^d)$ (\cite{Mi}). One may additionally assume that $\varphi\geq 0$. We apply the method presented in \cite{BGT3} and \cite{BGT6}. We start with the inequality

\begin{equation}
\label{eq:3.54}
P(|\langle\widetilde{X}_T,\varphi\otimes\psi\rangle|\geq \delta)\leq C\delta\int^{1/\delta}_0(1-{\rm Re}({ E\exp}\{-i\theta\langle\widetilde{X}_T,\varphi\otimes\psi\rangle\}))d\theta,
\end{equation}
valid for any $\psi\in {\cal S}(\erre),\delta>0$. Fix 
$0\leq t_1<t_2\leq 1$ and take $\psi$ approximating $\delta_{t_2}-\delta_{t_1}$ such that $\chi(t)=\int^1_t\psi(s)ds$ satisfies
\begin{equation}
\label{eq:3.55}
0\leq \chi\leq \UNO_{[t_1,t_2]}.
\end{equation}
Then the left hand side of (3.54) approximates
$$P(|\langle X_T(t_2),\varphi\rangle-\langle X_T(t_1),\varphi\rangle|\geq \delta).$$
So, in order to show tightness one should prove that the right hand side of 
(3.54) is estimated by
$$C(t^h_2-t^h_1)^{1+\sigma}\,\,\hbox{\rm for some}\,\, h,\sigma >0.$$
By the argument in \cite{BGT6} this reduces to showing that
\begin{equation}
\label{eq:3.56}
A(T)\leq C(t^h_2-t^h_1)^{1+\sigma}
\end{equation}
and
\begin{equation}
\label{eq:3.57}
I_1(T)\leq C(t^h_2-t^h_1)^{1+\sigma},
\end{equation}
where $I_1$ is defined by (\ref{eq:3.9}), and
\begin{equation}
\label{eq:3.58}
A(T)=\frac{H_T}{F^2_T}\int_{\errita^d}\int^T_0\int^s_0{\cal T}_{T-s}(\varphi{\cal T}_{s-u}\varphi)(x)
\chi\left(1-\frac{s}{T}\right)\chi\left(1-\frac{u}{T}\right)|x|^{-\gamma}dudsdx.
\end{equation}
The proofs of (\ref{eq:3.56}) and (\ref{eq:3.57}) are quite involved and lengthy,
 therefore, as an example we present only the argument for the case $\gamma<\alpha$, which, together with (\ref{eq:2.4}) implies
\begin{equation}
\label{eq:3.59}
d<\alpha\frac{1+\beta}{\beta}.
\end{equation}
We start with (\ref{eq:3.57}). By self-similarity of $p_s$ we have
\begin{equation}
\label{eq:3.60}
\int_{\errita^d}p_s(x-y)|x|^{-\gamma}dx\leq Cs^{-\gamma/\alpha},\quad y\in \erre^d.
\end{equation}
Using this, (\ref{eq:3.37}), (3.35) and (\ref{eq:3.55})
 we obtain
\begin{equation}
\label{eq:3.61}
I_1(T)\leq C(W_1(T)+W_2(T)),
\end{equation}
where
\begin{eqnarray}
\label{eq:3.62}
W_1(T)&=&\int_{\errita^d}\int^{t_1}_0s^{-\gamma/\alpha}\left(\int^{t_2}_{t_1}\int_{\errita^d}p_{u-s}(y-z)\widetilde{\varphi}_T(z)dzdu\right)^{1+\beta}dsdy,\\
\label{eq:3.63}
W_2(T)&=&\int_{\errita^d}\int^{t_2}_{t_1}s^{-\gamma/\alpha}\left(\int^{t_2}_s\int_{\errita^d}p_{u-s}(y-z)\widetilde{\varphi}_T(z)dzdu\right)^{1+\beta}dsdy.
\end{eqnarray}

Fix any $\rho$ such that
\begin{equation}
\label{eq:3.64}
\max\left\{\frac{d}{\alpha}-\frac{1}{\beta},0\right\}<\rho<1
\end{equation}
(see (\ref{eq:3.59}). For any fixed $s\in[0,t_1]$ we apply the Jensen inequality to the measure
$$\frac{(u-s)^{-\rho}}{\int^{t_2}_{t_1}(r-s)^{-\rho}dr}\UNO_{[t_1,t_2]}(u)du,$$
obtaining
$$W_1(T)\leq\int_{\errita^d}\int^{t_1}_0s^{-\gamma/\alpha}\left(\int^{t_2}_{t_1}(r-s)^{-\rho}dr\right)^\beta\int^{t_2}_{t_1}(u-s)^{-\rho}((u-s)^\rho p_{u-s}*\widetilde{\varphi}_T(y))^{1+\beta}dudsdy.$$
We have
\begin{equation}
\label{eq:3.65}
||p_{u-s}*\widetilde{\varphi}_T||_{1+\beta}\leq||p_{u-s}||_{1+\beta}||\varphi||_1=(u-s)^{-(d/\alpha)\beta/(1+\beta)}||p_1||_{1+\beta}||\varphi||_1,
\end{equation}
hence
$$W_1(T)\leq C(t_2-t_1)^{(1-\rho)\beta}\int^{t_2}_{t_1}\int^{t_1}_0s^{-\gamma/\alpha}(u-s)^{\rho\beta-(d/\alpha)\beta}dsdu,$$
which, after the substitution $s'=s/u$, by (\ref{eq:3.64}) and $\gamma<\alpha$, yields
\begin{equation}
\label{eq:3.66}
W_1(T)\leq C_1(t_2-t_1)^{(1-\rho)\beta}(t^h_2-t^h_1)\leq C_2(t^h_2-t^h_1)^{1+(1-\rho)\beta},
\end{equation}
where $h=2-\gamma/\alpha +\rho\beta-(d/\alpha)\beta>0$ by assumptions.

Next, by (\ref{eq:3.63}) we have
$$W_2(T)\leq\int_{\errita^d}\int^{t_2}_{t_1}
s^{-\gamma/\alpha}\left(\left(
\int^{t_2-t_1}_0p_udu\right)*\widetilde{\varphi}_T(y)\right)^{1+\beta}dsdy.$$
The Young inequality, substitution $u'=u/(t_2-t_1)$, and self-similarity imply
\begin{eqnarray}
W_2(T)&\leq&C( t^{1-\gamma/\alpha}_2-t^{1-\gamma/\alpha})
(t_2-t_1)^{1+\beta-(d/\alpha)\beta}||f||^{1+\beta}_{1+\beta}||
\varphi||^{1+\beta}_1\nonumber\\
\label{eq:3.67}
&\leq&C_1(t^{1-\gamma/\alpha}_2-t^{1-\gamma/\alpha}_1)^{2+\beta-(d/\alpha)\beta},
\end{eqnarray}
by (\ref{eq:3.41}).

Combining (\ref{eq:3.66}), (\ref{eq:3.67}), (\ref{eq:3.61}) and using 
(\ref{eq:3.59}),
 we obtain (\ref{eq:3.57}).

It remains to prove (\ref{eq:3.56}).

 Applying the usual substitutions to $A(T)$ given by (\ref{eq:3.58}) we obtain
\begin{eqnarray*}
\lefteqn{A(T)}\\
&=&\frac{H_TT^{2-d/\alpha-\gamma/\alpha}}{F^2_T}\int^1_0\int^1_s\int_{\errita^{3d}}|x|^{-\gamma}p_s(x-y)\widetilde{\varphi}_T(y)p_{u-s}(y-z)\widetilde{\varphi}_T(z)\chi(s)\chi(u)dxdzdyduds,
\end{eqnarray*}
hence, by (\ref{eq:3.60}), (\ref{eq:3.55}) and the H\"older inequality,
$$A(T)\leq C\frac{H_TT^{2-d/\alpha-\gamma/\alpha}}{F^2_T}
\int^{t_2}_{t_1}\int^{t_2}_ss^{-\gamma/\alpha}||\widetilde{\varphi}_T||_{\frac{1+\beta}{\beta}}||p_{u-s}*\widetilde{\varphi}_T||_{1+\beta}duds.$$
We use (\ref{eq:3.65}), (\ref{eq:3.45}) and (\ref{eq:2.5}), obtaining
$$A(T)\leq C_1(t^{1-\gamma/\alpha}_2-t^{1-1/\alpha}_1)(t_2-t_1)^{1-(d/\alpha)\beta/(1+\beta)},$$
which implies (\ref{eq:3.56}) by (\ref{eq:3.59}). This completes the proof of tightness. $\hfill\Box$
\vglue.5cm
\noindent
{\bf 3.4 Proof of Theorem 2.1(b)}
\vglue.5cm
We prove the theorem for
\begin{equation}
\label{eq:3.68}
\gamma=d<\frac{2+\beta}{1+\beta}\alpha
\end{equation}
(see Remark 2.2(c)). Recall that in this case $k(T)$ occurring in (\ref{eq:2.9})
and (\ref{eq:2.10}) is $\log T$.

According to the discussion in Section 3.1 it suffices to prove 
(\ref{eq:3.12}), (\ref{eq:3.13}) and (\ref{eq:3.15}). By the form of the limit process (see (\ref{eq:2.3})), (\ref{eq:3.12}) is equivalent to
\begin{equation}
\label{eq:3.69}
\lim_{T\to\infty}I_1(T)=\frac{\sigma(S_{d-1})}{\alpha}\int_{\errita^d}\int^1_0p_s(y)
\biggl(\int^1_sp_{u-s}(y)\chi(u)du\biggr)^{1+\beta}dsdy\left(\int_{\errita^d}\varphi(z)dz\right)^{1+\beta},
\end{equation}
where $\sigma(S_{d-1})$ is the measure of the unit sphere in $\erre^d (=2$ if $d=1)$.

By (\ref{eq:3.9}), (\ref{eq:3.2}), (\ref{eq:3.6}),
 using similar substitutions as in the previous section, we obtain
\begin{eqnarray}
\lefteqn{I_1(T)} \nonumber\\
\label{eq:3.70}
&=&\frac{1}{\log T}\int_{\errita^d}\int^1_0\int_{\errita^d}p_s(x
T^{-1/\alpha}-y)\biggl(\int^1_s\int_{\errita^d}p_{u-s}(y-z)\chi(u)\widetilde{\varphi}_T(z)dzdu\biggr)^{1+\beta}\frac{1}{1+|x|^d}dydsdx,\nonumber \\
\end{eqnarray}
where $\widetilde{\varphi}_T$ is given by (\ref{eq:3.34}). We write
\begin{equation}
\label{eq:3.71}
I_1(T)=I'_1(T)+I^{''}_1(T)+I_1^{'''}(T),
\end{equation}
where
\begin{eqnarray}
\label{eq:3.72}
I'_1(T)&=&\frac{1}{\log T}\int_{1<|x|<T^{1/\alpha}}\int^1_0\int_{\errita^d} ...
\\
\label{eq:3.73}
I^{''}_1(T)&=&\frac{1}{\log T}\int_{|x|\geq T^{1/\alpha}}\int^1_0\int_{\errita^d} ...\\
\label{eq:3.74}
I^{'''}_1(T)&=&\frac{1}{\log T}\int_{|x|\leq 1}\int^1_0\int_{\errita^d} ... .
\end{eqnarray}

Passing to polar coordinates in the integral with respect to $x$ we have
$$I'_1(T)=\frac{1}{\log T}\int^{T^{1/\alpha}}_1\int_{S_{d-1}}\int^1_0\int_{\errita^d}p_s(wrT^{-1/\alpha}-y)(g_s*\widetilde{\varphi}_T(y))^{1+\beta}\frac{r^{d-1}}{1+r^d}dyds\sigma(dw)dr,$$
where $g$ is defined by (\ref{eq:3.35}). The crucial step is the substitution
\begin{equation}
\label{eq:3.75}
r'=\frac{\log r}{\log T},
\end{equation}
which gives
$$I'_1(T)=\int^{1/\alpha}_0\int_{S_{d-1}}\int^1_0\int_{\errita^d}p_s(wT^{r-1/\alpha}-y)(g_s*\widetilde{\varphi}_T(y))^{1+\beta}\frac{T^{rd}}{1+T^{rd}}dyds\sigma(dw)dr.$$
It is now clear that if one could pass to the limit under the integrals as $T\to\infty$, then $I'_1(T)$ would converge to the right hand side of (\ref{eq:3.69}). This procedure is indeed justified by the fact that for $f$ defined by 
(\ref{eq:3.22}) we have
\begin{equation}
\label{eq:3.76}
f\in L^{2+\beta}(\erre^d),
\end{equation}
which follows from
 (\ref{eq:3.68}), (\ref{eq:3.23}) and (\ref{eq:3.24}). We omit the details, which are similar to the argument in \cite{BGT6} (see (\ref{eq:3.51}) therein).

Next we show that $I^{''}_1(T) $ and $I^{'''}_1(T)$ tend to zero. In $I^{''}_1(T)$ (see (\ref{eq:3.73})) we substitute $x'=xT^{-1/\alpha}$ and we use
(\ref{eq:3.36}), obtaining
\begin{eqnarray*}
I^{''}_1(T)&\leq&\frac{C}{\log T}\int_{|x|>1}\int^1_0\int_{\errita^d}p_s(x-y)(f*\widetilde{\varphi}_T(y))^{1+\beta}\frac{T^{d/\alpha}}{1+|x|^dT^{d/\alpha}}dydsdx\\
&\leq&\frac{C_1}{\log T}||f*\widetilde{\varphi}_T||^{1+\beta}_{1+\beta}
\,\,\leq\,\,\frac{C_1}{\log T}||f||^{1+\beta}_{1+\beta}||\varphi||^{1+\beta}_1\to 0.
\end{eqnarray*}

$I^{'''}_1(T)$ (see (\ref{eq:3.74})) is estimated as follows:
\begin{eqnarray*}
I^{'''}_1(T)&\leq&\frac{C}{\log T}\int_{|x|\leq 1}f*(f*\widetilde{\varphi}_T)^{1+\beta}(xT^{-1/\alpha})dx\\
&\leq&\frac{C_2}{\log T}||f||_{2+\beta}||(f*\widetilde{\varphi}_T)^{1+\beta}||_{(2+\beta)/(1+\beta)}\\
&\leq&\frac{C_2}{\log T}||f||^{2+\beta}_{2+\beta}||\varphi||^{1+\beta}_1\to 0,
\end{eqnarray*}
by (\ref{eq:3.76}). This and (\ref{eq:3.71}) prove (\ref{eq:3.69}).

To prove (\ref{eq:3.13}) we use (\ref{eq:3.14}) and easily obtain
$$I_2(T)\leq\frac{CH_TT^{2-2d/\alpha}}{F^2_T}\int_{\errita^2}f*(\widetilde{\varphi}_T(f*\widetilde{\varphi}_T))(xT^{-1/\alpha})\frac{1}{1+|x|^d}dx.$$
We write the right-hand side as the sum of integrals over $\{|x|\leq T^{1/\alpha}\}$ and $\{|x|>T^{1/\alpha}\}$. To estimate the integral over $\{|x|\leq T^{1/\alpha}\}$ we use
\begin{equation}
\label{eq:3.77}
\sup_{T>2}\frac{1}{\log T}
\int_{|x|\leq T^{1/\alpha}}\frac{1}{1+|x|^d}dx<\infty,
\end{equation}
and in the second integral we apply $1/(1+|x|^d)\leq T^{-d/\alpha}$.
For each of the integrals we use appropriately the H\"older inequality, properties of the convolution and (\ref{eq:3.45}), obtaining the estimates
 $C_1T^{2((d/\alpha)1/(2+\beta)
-1/(1+\beta))}$ and $C_2T^{(d/\alpha)1/(2+\beta)-2/(1+\beta)}$, respectively (the factors involving negative powers of $H_T$ and $\log T$ have been estimated by constants). These bounds tend to zero as $T\to\infty$ by 
(\ref{eq:3.68}). We omit details. This proves (\ref{eq:3.13}).

To prove (\ref{eq:3.16}) and (\ref{eq:3.17}) we use (\ref{eq:3.18}) and
(\ref{eq:3.19}). Again, we consider separately the integrals over $\{|x|\leq T^{1/\alpha}\}$ and $\{|x|>T^{1/\alpha}\}$, and apply the same tricks as for $I_2(T)$.

For $J_1(T)$ we obtain the estimate
$$J_1(T)\leq CT^{(d/\alpha)(1+\beta)/(2+\beta)-1}\to 0$$
($\log T$ and $H_T$ appear with negative powers only), whereas
$$J_2(T)\leq C_1\frac{T}{H^\beta_T(\log T)^\beta}+C_2\frac{T}{H^\beta_T(\log T)^{1+\beta}}\to 0$$
by assumption (\ref{eq:2.10}). The proof of Theorem 2.1 is complete $\hfill\Box$\vglue.5cm
\noindent
{\bf 3.5 Proof of Proposition 2.3}
\vglue.5cm
Properties (a)-(c) are clear, following from (\ref{eq:2.1}) and Theorem 2.1(a). Recall that the index of self-similarity is defined as $a\in\erre$ such that the process $(\xi_{ct})_{t\in\errita_+}$ has the same distribution as $(c^a\xi_t)_{t\in \errita_+}$ for any $c>0$.

To calculate the dependence exponent of $\xi$ (see (\ref{eq:2.11}), (\ref{eq:2.12})) first note that by (\ref{eq:2.1}) and Proposition 3.4.2 of \cite{ST} the finite-dimensional distributions of $\xi$ are given by 
\begin{eqnarray}
\lefteqn{E{\rm exp}\{i(z_1\xi_{t_1}+\cdots+z_k\xi_{t_k})\}}\nonumber\\
&=&{\rm exp}\left\{-\int_{\errita^{d+1}}
\right.
\left[\biggl|\sum^k_{j=1}z_j\left(\int_{\errita^d}
p_r(x-y)|y|^{-\gamma}dy\right)^{1/(1+\beta)}
\UNO_{[0,t_j]}(r)\int^{t_j}_rp_{u-r}(x)du\biggr|^{1+\beta}\right.\nonumber\\
\label{eq:3.78}
&&\times \Biggl(1-i\,\,{\rm sgn}\Biggl(\sum^k_{j=1}z_j\left(
\int_{\errita^d}p_r(x-y)|y|^{-\gamma}dy\right)^{1/(1+\beta)}
\bigg. \biggr.\nonumber\\
&&\times \left.\left.\left.
\UNO_{[0,t_j]}(r)\int^{t_j}_rp_{u-r}(x)du\right){\rm tan}
\frac{\pi}{2}(1+\beta)\Bigg) \right]drdx\right\}.
\end{eqnarray}

The argument goes along the lines of the proof of Theorem 2.7 of \cite{BGT3}. For fixed $z>0$ and $0\leq u<v<s<t$ we define $D^+_T=D_T(1,z;u,v,s,t)$ and $D^-_T=D_T(1,-z;u,v,s,t)$ (the formulas for $D^+,D^-$ are obtained from (\ref{eq:2.12}) and (\ref{eq:3.78})), and we prove
$$
D^\pm_T\leq\left\{
\begin{array}{l}
CT^{-d/\alpha}\,\,\hbox{\rm if either}\quad \alpha=2\,\,{\rm or}\,\,\beta>
{(d-\gamma)}/{(d+\alpha)},\\
CT^{-(d/\alpha)\delta}\,\,\hbox{\rm for any}\quad \beta<\delta<1+\beta+
{(d-\gamma)}/{(d+\alpha)}\,\,{\rm if}\,\,\alpha<2,\beta\leq
{(d-\gamma)}/{(d+\alpha)},
\end{array}\right.
$$
and for $T$ sufficiently large,
\begin{eqnarray}
D^+_T&\geq&CT^{-d/\alpha},\phantom{333ddsdfasdfasdfasdfsdfsdddddddddddddddddddddddddddddddddddd33}\nonumber\\
\label{eq:3.79}
D^+_T&\geq&CT^{-(d/\alpha)\delta}\,\,\hbox{\rm for any}\,\,\delta>1+\beta-
{(d-\gamma)}/{(d+\alpha)}\,{\rm if}\,\,\alpha<2,\beta\leq
({d-\gamma})/({d+\alpha}).
\end{eqnarray}

The upper estimates are obtained similarly as 
(4.3), (4.4) in \cite{BGT3} and (\ref{eq:3.107}) in \cite{BGT6}. The only difference is that in formulas (4.9) and (4.10) in \cite{BGT3} a new factor, $\int_{\errita^d}p_r(x-y)|y|^{-\gamma}dy$, appears (which 
 corresponds to $p_r(x)$ in (\ref{eq:3.109}) in \cite{BGT6}). This factor is responsible for the new long-range dependence threshold and the form of the dependence exponent (\ref{eq:2.13}). In the estimates we use (\ref{eq:3.30}).

The first of the lower estimates is obtained exactly as (4.18) in \cite{BGT3}. The new expression,
$$\int^{(u+v)/2}_u\int_{|x|\leq 1}\int_{\errita^d}p_r(x-y)|y|^{-\gamma}dydxdr,$$
that appears at the right-hand side is finite by (\ref{eq:3.31}).

To derive (\ref{eq:3.79}) we argue as in (4.22),
(4.24) of \cite{BGT3} and we apply estimates (4.21) (which holds for $|x|\leq T^{1/\alpha}$) and (4.23) therein,  obtaining
\begin{equation}
\label{eq:3.80}
D^+_T\geq CT^{-(d/\alpha)(1+\beta)+\varepsilon\beta(d+\alpha)}
\int^{(u+v)/2}_{u+(v-u)/4}\int_{1\leq|x|\leq T^{d/(d+\alpha)\alpha-\varepsilon}}\int_{\errita^d}p_r(x-y)|y|^{-\gamma}dydxdr,
\end{equation}
where $\varepsilon>0$ is sufficiently small.

For $1\leq |x|\leq T^{d/(d+\alpha)\alpha-\varepsilon}$ we have
\begin{eqnarray*}
\int_{\errita^d}p_r(x-y)|y|^{-\gamma}dy&\geq& C|x|^{-\gamma}\int_{|x-y|\leq1/2}p_r(x-y)dy\\
&\geq&C_1|x|^{-\gamma}\inf_{\frac{v-u}{4}+u<r<\frac{u+v}{2}}
\inf_{|z|\leq\frac{1}{2}}p_r(z)\geq C_2T^{-d\gamma/d(d+\alpha)\alpha+\varepsilon\gamma}.
\end{eqnarray*}
Putting this into (\ref{eq:3.80}) we obtain (\ref{eq:3.79}). $\hfill\Box$
\vglue.5cm
\noindent
{\bf 3.6 Proof of Theorem 2.5}
\vglue.5cm
Each of the cases requires a different proof, and none of them is straightforward. We will present a detailed proof of the part (a) only. In the remaining cases we will confine ourselves to explaining why the limit processes have the forms given in the theorem (recall that part (b)(iv) has been proved in \cite{BGT6}).
It seems instructive to compare the proofs for this theorem to the argument given in the proof of Theorem 2.1(b) for $\gamma=d$. Although the critical cases are of different kinds, some of the technical tricks repeat in all cases, nevertheless they are applied in a slightly different way and are far from being identical.
\vglue.5cm
\noindent
{\bf Proof of case (a)} To simplify calculations we again consider the measure $\mu_\gamma$ of the form $\mu_\gamma(dx)=|x|^{-\gamma}dx$ instead of 
(\ref{eq:1.8}).

In (\ref{eq:3.9}) we substitute $u'=s-u$ and then $s'=(T-s)/T$, obtaining
$$
I_1(T)=\frac{H_TT}{F^{1+\beta}_T}\int_{\errita^d}\int^1_0\int_{\errita^d}p_{sT}(x-y)\left(\int^{T(1-s)}_0\int_{\errita^d}p_u(y-z)\varphi(z)
\chi\biggl(s+\frac{u}{T}\biggr)dzdu\right)^{1+\beta}|x|^{-\gamma}dydsdx.$$
Using (\ref{eq:2.15}), (\ref{eq:3.20}) and substitution $x'=x(sT)^{-1/\alpha}$ we have\\
\vbox{
\begin{eqnarray}
I_1(T)&=&\frac{1}{\log T}\int_{\errita^d}\int^1_0\int_{\errita^d}p_1(x-y(sT)^{-1/\alpha})\nonumber\\
&&\times\left(\int^{T(1-s)}_0\int_{\errita^d}p_u(y-z)\varphi(z)\chi
\biggl(s+\frac{u}{T}\biggr)dzdu\right)^{1+\beta}s^{-\gamma/\alpha}|x|^{-\gamma}dydsdx\nonumber\\
\label{eq:3.81}
&=&I'_1(T)+I^{''}_1(T)+I^{'''}_1(T),
\end{eqnarray}
}
where
\begin{eqnarray}
\label{eq:3.82}
I'_1(T)&=&\frac{1}{\log T}\int_{\errita^d}\int^1_0\int_{1\leq|y|\leq T^{1/\alpha}}\ldots\\
\label{eq:3.83}
I''_1(T)&=&\frac{1}{\log T}\int_{\errita^d}\int^1_0
\int_{|y|>T^{1/\alpha}}\ldots\\
\label{eq:3.84}
I'''_1(T)&=&\frac{1}{\log T}\int_{\errita^d}\int^1_0
\int_{|y|<1}\ldots
\end{eqnarray}
Passing to polar coordinates in the integral with respect to $y$ and making substitution (\ref{eq:3.75}) we obtain
\begin{eqnarray*}
\lefteqn{I'_1(T)}\\
&=&\int_{\errita^d}\int^1_0\int^{1/\alpha}_0\int_{S_{d-1}}p_1(x-wT^{r-1/\alpha}s^{-1/\alpha})\nonumber \\
&&\times\left(\int^{T(1-s)}_0\int_{\errita^d}p_u(wT^r-z)\varphi(z)\chi
\biggl(\frac{u}{T}+s\biggr)dzdu\right)^{1+\beta}
 s^{-\gamma/\alpha}|x|^{-\gamma}T^{rd}\sigma(dw)drdsdx.
\end{eqnarray*}
We substitute $z'=T^{-r}z, u'=uT^{-r\alpha}$, use (\ref{eq:3.20}) and 
(\ref{eq:2.14}), arriving at
\begin{equation}
\label{eq:3.85}
I'_1(T)=\int_{\errita^d}\int^1_0\int^{1/\alpha}_0\int_{S_{d-1}}p_1(x-ws^{-1/\alpha}T^{r-1/\alpha})\left(h_T(r,s,w)\right)^{1+\beta}s^{-\gamma/\alpha}|x|^{-\gamma}\sigma(dw)drdsdx,
\end{equation}
where
\begin{equation}
\label{eq:3.86}
h_T(r,s,w)=\int^{T^{1-r\alpha}(1-s)}_0\int_{\errita^d}p_u(w-z)T^{rd}\varphi(zT^r)\chi(s+uT^{r\alpha-1})dzdu.
\end{equation}

It is clear that on the set of integration one should have 
\begin{eqnarray}
\lim_{T\to\infty}h_T(r,s,w)&=&\int_{\errita^d}\varphi(z)dz\int^\infty_0p_u(w)du\chi(s)\nonumber\\
\label{eq:3.87}
&=&C_{d,\alpha}\int_{\errita^d}\varphi(z)dz\chi(s),
\end{eqnarray}
where $C_{d,\alpha}$ is given by (\ref{eq:1.13}), which should yield
\begin{equation}
\label{eq:3.88}
\lim_{T\to\infty}I'_1(T)=C^{1+\beta}_{d,\alpha}\frac{1}{\alpha}\sigma(S_{d-1})\int_{\errita^d}p_1(x)|x|^{-\gamma}dx\int^1_0s^{-\gamma/\alpha}\chi^{1+\beta}(s)ds\biggl(\int_{\errita^d}\varphi(z)dz\biggr)^{1+\beta}.
\end{equation}
However, (\ref{eq:3.87}) and (\ref{eq:3.88}) need a justification. It is easy to see that the first integral in (\ref{eq:3.86}) can be replaced by $\int^\infty_0 du$. Since
$$\lim_{T\to\infty}\int_{\errita^d}p_u(w-z)T^{rd}\varphi(zT^r)dz=p_u(w)\int_{\errita^d}\varphi(z)dz,$$
it is clear that in order to prove (\ref{eq:3.87}) it suffices to show that
$$\sup_{T>2}\sup_{w\in S_{d-1}}\int_{\errita^d}p_u(w-z)T^{rd}\varphi(zT^r)dz$$
is integrable in $u$. This
 is clear for $u\geq 1$ since $d>\alpha$, and for $u<1$ we argue similarly as in (\ref{eq:3.39}) obtaining an integrable bound $C_1p_u(w/2)+C_2$.
 In the same way one shows that $h_T(r,s,w)\leq C$. This together with (\ref{eq:3.87}) easily implies (\ref{eq:3.88}).

Next, it is easy to see that for $I'''_1$ defined by (\ref{eq:3.84}) we have
\begin{equation}
\label{eq:3.89}
I'''_1(T)\leq\frac{C}{\log T}\to 0.
\end{equation}
A little more work is needed to prove that also
\begin{equation}
\label{eq:3.90}
\lim_{T\to\infty}I''_1(T)=0.
\end{equation}
By (\ref{eq:3.83}),
\begin{eqnarray*}
I''_1(T)&\leq&\frac{C}{\log T}\int_{|y|>T^{1/\alpha}}\left(\int^T_0\int_{\errita^d}p_u(y-z)\varphi(z)dzdu\right)^{1+\beta}dy\\
&\leq&C_1(R_1(T)+R_2(T)),
\end{eqnarray*}
where
\begin{eqnarray*}
R_1(T)&=&\frac{1}{\log T}\int_{|y|>T^{1/\alpha}}\left(\int^T_0\int_{|z|\leq\frac{T^{1/\alpha}}{2}}p_u(y-z)\varphi(z)dzdu\right)^{1+\beta}dy,\\
R_2(T)&=&\frac{1}{\log T}\int_{|y|>T^{1/\alpha}}
\left(\int^T_0\int_{|z|>\frac{T^{1/\alpha}}{2}}
p_u(y-z)\varphi(z)dzdu\right)^{1+\beta}dy.
\end{eqnarray*}
We have
\begin{eqnarray*}
R_1(T)&\leq&\frac{1}{\log T}\int_{|y|>T^{1/\alpha}}\left(\int^T_0p_u
\left(\frac{y}{2}\right)du\right)^{1+\beta}dy\left(\int_{\errita^d}\varphi(z)dz\right)^{1+\beta}\\
&=&\frac{C}{\log T}\int_{|y|>1}\left(\int^1_0p_u\left(\frac{y}{2}\right)du\right)^{1+\beta}dy,
\end{eqnarray*}
after obvious substitutions and using (\ref{eq:2.14}). Hence $\lim\limits_{T\to\infty}R_1(T)=0$ by (\ref{eq:3.23}). Furthermore,
\begin{eqnarray*}
R_2(T)&\leq&\frac{C}{\log T}\int_{|y|>T^{1/\alpha}}
\left(\int^T_0\int_{|z|>\frac{T^{1/\alpha}}{2}}p_u(y-z)\frac{|z|^2\varphi(z)}{T^{2/\alpha}}dzdu\right)^{1+\beta}dy\\
&\leq&\frac{C_1}{T^{(2/\alpha)(1+\beta)}},
\end{eqnarray*}
since, under (\ref{eq:2.14}),
$$\sup_{T}\frac 1 {\log T} \int_{\errita^d}\left(\int^T_0\int_{\errita^d}p_u(y-z)\varphi_1(z)dzdu\right)^{1+\beta}dy<\infty,$$
for any $\varphi_1\in{\cal S}(\erre^d)$ by (\ref{eq:3.33}) in \cite{BGT4} (in our case $\varphi_1(z)=|z|^2\varphi(z))$. This proves (\ref{eq:3.90}), and by
(\ref{eq:3.81})-(\ref{eq:3.84}), (\ref{eq:3.88}) and (\ref{eq:3.89}) we have established (\ref{eq:3.12}).

To prove (\ref{eq:3.13}) we use (\ref{eq:3.14}), which, after standard 
substitutions, gives 
\begin{eqnarray*}
\lefteqn{I_2(T)}\\
&\leq&CH^{-1/(1+\beta)}_T
T^{-2/(1+\beta)+2\gamma/\alpha(1+\beta)+2-d/\alpha-\gamma/\alpha}(\log T)^{-2/(1+\beta)} \int_{\errita^d}f_\gamma(y)\widetilde{\varphi}_T(y)
(f*\widetilde{\varphi}_T)(y)dy
\end{eqnarray*}
(see (\ref{eq:2.15}), (\ref{eq:3.22}), (\ref{eq:3.29}), (\ref{eq:3.34})). By 
(\ref{eq:3.31}) we have
\begin{equation}
\label{eq:3.91}
I_2(T)\leq C_1T^{-2/(1+\beta)+2\gamma/\alpha(1+\beta)+2-d/\alpha-\gamma/\alpha}
\sup_y f*\widetilde{\varphi}_T(y).
\end{equation}
The assumptions (\ref{eq:2.14}) and $\gamma<\alpha$ imply that
$$\frac{2}{1+\beta}-\frac{2\gamma}{\alpha(1+\beta)}-2+\frac{d}{\alpha}+\frac{\gamma}{\alpha}=\frac{1}{\beta}\theta$$
for some $\theta>1$. Observe that
\begin{equation}
\label{eq:3.92}
f\in L^q(\erre^d),\,\, 1\leq q<1+\beta,
\end{equation}
by (\ref{eq:3.23}), (\ref{eq:3.24}) and (\ref{eq:2.14}). Fix $q>1$ such that $(1+\beta)/\theta<q<1+\beta$ and $p=q/(q-1)$. By the H\"older inequality and 
(\ref{eq:3.45}) we obtain
$$I_2(T)\leq C_1T^{-\theta/\beta+((1+\beta)/\beta)(1/q)}||\varphi||_p\to 0$$
as $T\to\infty$, by assumption on $q$.

To prove (\ref{eq:3.16}) we use (\ref{eq:3.18}) which, by an analogous argument as for $I_2$, gives
\begin{equation}
\label{eq:3.93}
J_1(T)\leq CT^{-1+\gamma/\alpha-(1+\beta)/\beta}||f*(\widetilde{\varphi}_T(f*\widetilde{\varphi}_T))||^{1+\beta}_{1+\beta}.
\end{equation}
Taking $q<1+\beta$ sufficiently close to $1+\beta$ and using (\ref{eq:3.92}), the Young and H\"older inequalities can be applied to conclude that, by 
(\ref{eq:3.45}),
$$||f*(\widetilde{\varphi}_T(f*\widetilde{\varphi}_T))||^{1+\beta}_{1+\beta}\leq O(T^r),$$
where $r<1-\gamma/\alpha+(1+\beta)/\beta$; we omit details. This and 
(\ref{eq:3.93}) yield (\ref{eq:3.16}).

To prove (\ref{eq:3.17}) we write an estimate similar to (\ref{eq:3.50}), 
namely,
\begin{equation}
\label{eq:3.94}
J_2(T)\leq CT^{\beta(\gamma/\alpha-1)}R(T),
\end{equation}
where $R(T)$ is defined by (\ref{eq:3.51}). In the present case (\ref{eq:3.52})
does not hold, but similarly as before, using the Young inequality, 
(\ref{eq:3.92}) and (\ref{eq:3.45}) it can be shown that $R(T)=O(T^\varepsilon)$ for any $\varepsilon >0$. This and (\ref{eq:3.43}) imply (\ref{eq:3.17}) since $\gamma<\alpha$. The proof of part (a) of the theorem is complete.
\vglue.5cm
\noindent
{\bf Sketch of the proof of case (b)}
\vglue.5cm
(i) We repeat the argument as in (\ref{eq:3.81})
-(\ref{eq:3.86}) with $F_T$ given by (\ref{eq:2.16}). Again, it can be shown that the only significant term is $I'_1(T)$, i.e., $\lim_{T\to\infty}I_1(T)=\lim_{T\to\infty}I'_1(T)$. In order to derive this limit, in (\ref{eq:3.85}) we substitute $s'=sT^{1-r\alpha}$, obtaining
$$I'_1(T)=\frac{1}{\log T}\int_{\errita^d}\int^{1/\alpha}_0\int^{T^{1-r\alpha}}_0\int_{S_{d-1}}p_1(x-ws^{1/\alpha})(h_T(r,sT^{-1+r\alpha},w))^{1+\beta}s^{-1}|x|^{-\alpha}\sigma(dw)dsdrdx.$$
It is easy to see that the limit remains the same if the integral 
$\int^{T^{1-r\alpha}}_0\ldots ds$ is replaced by $\int^{T^{1-r\alpha}}_1\ldots ds$. Let $\widetilde{I}'_1(T)$ denote $I'_1(T)$ after this change. Next, we substitute $s'=\log s/\log T$ and we have
$$
\widetilde{I}'_1(T)=
\int^{1/\alpha}_0\int_{S_{d-1}}\int^{1-r\alpha}_0\int_{\errita^d}
p_1(x-wT^{-s/\alpha})|x|^{-\alpha}(h_T(r,T^{s-1+r\alpha},
w))^{1+\beta}dxds\sigma(dw)dr.
$$
By (\ref{eq:3.86}), it is clear that on the set of integration 
$$\lim_{T\to\infty} h_T(r,T^{s-1+r\alpha},w)=C_{d,\alpha}\int_{\errita^d}\varphi(z)dr\chi(0).$$
This shows that we should have
\begin{eqnarray}
\lefteqn{\lim_{T\to\infty}I_1(T)=
\lim_{T\to\infty}\widetilde{I}'_1(T)}\nonumber\\
\label{eq:3.95}
&=&C^{1+\beta}_{d,\alpha}\sigma(S_{d-1})\int_{\errita^d}p_1(x)|x|^{-d}dx\int^{1/\alpha}_0(1-r\alpha)dr
\biggl(\int_{\errita^d}\varphi(z)dz\biggr)^{1+\beta}\chi^{1+\beta}(0),
\end{eqnarray}
and this passage to the limit can be indeed justified. The right-hand side of 
(\ref{eq:3.95}) is equal to 
$\log E {\rm exp}{\{-C\langle\widetilde{X},\varphi\otimes\psi\rangle\}},$
where $X(=K\lambda\vartheta)$ is the limit process defined in the theorem. We skip the remaining parts of the proof.

(ii) In 
(\ref{eq:3.9}) we substitute $u'=s-u$ and then $s'=T-s$, obtaining
$$I_1(T)=I'_1(T)+I''_1(T),$$
where
\begin{eqnarray*}
I'_1(T)
&=&\frac{1}{\log T}\int^T_1\int_{\errita^d}\int_{\errita^d}p_1((x-y)s^{-1/\alpha})s^{-d/\alpha}|x|^{-\gamma}\\
&\times&\left(\int^{T-s}_0\int_{\errita^d}p_u(y-z)\varphi(z)\chi
\biggl(\frac{u+s}{T}\biggr)dzdu\right)^{1+\beta}dydxds
\end{eqnarray*}
and
$$I''_1(T)=\frac{1}{\log T}\int^1_0\int_{\errita^d}\int_{\errita^d}\ldots dydxds.$$
It can be shown that
$$\lim_{T\to\infty}I''_1(T)=0.$$
In $I'_1(T)$ we substitute $x'=xs^{-1/\alpha},y'=ys^{-1/\alpha},u'=
{u}/{s}$, and use (\ref{eq:3.20}), which gives
\begin{eqnarray*}
I'_1(T)&=&\frac{1}{\log T}\int^T_1\int_{\errita^d}\int_{\errita^d}p_1(x-y)|x|^{-\gamma}s^{-(d/\alpha)\beta-\gamma/\alpha+(1+\beta)}\\
&&\times\biggl(\int^{T/s-1}_0\int_{\errita^d}p_u(y-zs^{-1/\alpha})\varphi(z)
\chi\biggl(\frac{s(u+1)}{T}\biggr)dzdu\biggr)^{1+\beta}dydxds.
\end{eqnarray*}
By (\ref{eq:2.18}), $s^{-(d/\alpha)\beta-\gamma/\alpha+1+\beta}=s^{-1}$, so, the substitution $s'=\log s/\log T$ yields
\begin{eqnarray}
I'_1(T)&=&\int^1_0\int_{\errita^d}\int_{\errita^d}p_1(x-y)
|x|^{-\gamma}\nonumber\\
\label{eq:3.96}
&&\times\left(\int^{T^{1-s}-1}_0\int_{\errita^d}p_u(y-zT^{-s/\alpha})\varphi(z)
\chi\left(((u+1)T^{s-1}\right)dzdu\right)^{1+\beta}dydxds.\nonumber\\
\end{eqnarray}
It is now seen that one should have
\begin{eqnarray}
\lefteqn{\lim_{T\to\infty}I_1(T)=\lim_{T\to\infty}I'_1(T)}\nonumber\\
\label{eq:3.97}
&=&C^{1+\beta}_{\alpha,d}\int_{\errita^d}\int_{\errita^d}p_1(x-y)|x|^{-\gamma}|y|^{-(d-\alpha)(1+\beta)}dxdy
\left(\int_{\errita^d}\varphi(z)dz\right)^{1+\beta}\chi^{1+\beta}(0)
\end{eqnarray}
Note that the integrals are finite by (\ref{eq:3.21}), (\ref{eq:3.28}) 
and 
(\ref{eq:2.18}). The justification of (\ref{eq:3.97})  requires some work, but we omit it for brevity.

(iii) As $d=\gamma$, we must keep the measure $\mu_\gamma$ in its original form (\ref{eq:1.8}). 

Arguing as in the proof of (ii) and taking into account 
(\ref{eq:2.16}), instead 
of (\ref{eq:3.96}) we obtain
\begin{eqnarray*}
\lefteqn{I'_1(T)
=\frac{1}{\log T}\int^1_0\int_{\errita^d}\int_{\errita^d}p_1(x-y)
\frac{T^{sd/\alpha}}{1+|xT^{s/\alpha}|^d}}\\
&&\times\left(\int^{T^{1-s}-1}_0\int_{\errita^d}p_u(y-z T^{-s/\alpha})\varphi(z)
\chi((u+1)T^{s-1})dzdu\right)^{1+\beta}dydxds.
\end{eqnarray*}

Since
$$\lim_{T\to\infty}\frac{1}{\log T}\int_{\errita^d}p_1(x-y)
\frac{T^{sd/\alpha}}{1+|xT^{s/\alpha}|^d}dx=s\frac{1}{\alpha}\sigma(S_{d-1})p_1(y),$$
it can be shown, with some effort, that
\begin{eqnarray*}
\lim_{T\to\infty}I_1(T)&=&\lim_{T\to\infty}I'_1(T)\\
&=&C^{1+\beta}_{\alpha,d}\frac{1}{2\alpha}\sigma(s_{d-1})\int_{\errita^d}p_1(y)|y|^{-(d-\alpha)(1+\beta)}dy\biggl(\int_{\errita^d}\varphi(z)dz\biggr)^{1+\beta}\chi^{1+\beta}(0).
\end{eqnarray*}
Again, we omit the remaining parts of the proof.
\vglue.5cm
\noindent
{\bf 3.7 Proof of Theorem 2.6}
\vglue.5cm
We only give an 
outline of the proof. The following lemma is constantly used.

\vglue.5cm
\noindent
{\bf Lemma} {\it 
Let $\varphi\in{\cal S}(\erre^d),\varphi\geq 0$.

\noindent
(a) If $d>\alpha(2+\beta)/(1+\beta)$, then the functions $G\varphi, G(G\varphi)^{1+\beta}$ and $G(G\varphi)^2$ are bounded.

\noindent
(b) If $d>\alpha{(1+\beta)}/{\beta}$, then additionally 
$(G\varphi)^{1+\beta}$ and 
$(G(G\varphi)^{1+\beta})^{1+\beta}$ are integrable (and bounded).

\noindent
(c) If $\alpha<\gamma\leq d$ and $d>\alpha(2+\beta)/\beta-\gamma/\beta$, then
 additionally to the properties in (a),}
$$\int_{\errita^d}G(G\varphi)^{1+\beta}(x)\frac{1}{1+|x|^\gamma}dx<\infty.$$

This Lemma follows easily from (\ref{eq:1.12}) and (\ref{eq:3.25})-(\ref{eq:3.28}). 
\vglue.5cm
\noindent
{\bf Proof of part (a) of the theorem.} As before we consider $\mu_\gamma(dx)=|x|^{-\gamma}dx$. In (\ref{eq:3.9}) we substitute $u'=s-u$, then $s'=(T-s)/T$ and, finally, $x'=xT^{1/\alpha}s^{-1/\alpha}$, obtaining
\begin{eqnarray*}
\lefteqn{I_1(T)}\\
&=&\int_{\errita^d}\int^1_0\int_{\errita^d}p_1(x-ys^{-1/\alpha}T^{-1/\alpha})\left(\int^{T(1-s)}_0{\cal T}_u\varphi(y)\chi\biggl(s+\frac{u}{T}\biggr)du\right)^{1+\beta}s^{-\gamma/\alpha}|x|^{-\gamma}dydsdx
\end{eqnarray*}
(see (\ref{eq:2.24})). It is easily seen that by part (b) of the Lemma we have
\begin{equation}
\label{eq:3.98}
\lim_{T\to\infty}I_1(T)=\int_{\errita^d}p_1(x)|x|^{-\gamma}dx\int^1_0s^{-\gamma/\alpha}\chi^{1+\beta}(s)ds\int_{\errita^d}(G\varphi(y))^{1+\beta}dy.
\end{equation}
For $\beta<1$ this is exactly $\log E{\rm exp}\{-C\langle\widetilde{X},\varphi\otimes\psi\rangle\}$, where $X$ is the limit process described in the theorem.
Moreover, in this case (\ref{eq:3.14}) and boundedness of $G\varphi$ easily imply
$$I_2(T)\leq C \quad T^{(1-\gamma/\alpha)(1-2/(1+\beta))}\to 0.$$
For $\beta=1$ we use (\ref{eq:3.10}) and (\ref{eq:3.4}), obtaining
$$I_2(T)=I'_2(T)-I''_2(T)-\frac{V}{2}I'''_2(T),$$
where
\begin{eqnarray*}
I'_2(T)&=&\frac{H_T}{F^2_T}\int^T_0\int_{\errita^d}\int_{\errita^d}p_s(x-y)|x|^{-\gamma}dx\varphi(y)\chi(s)\int^T_s{\cal T}_{u-s}\varphi(y)\chi
\left(\frac{u}{T}\right)dudyds,\\
I''_2(T)&=&\frac{H_T}{F^2_T}\int_{\errita^d}\int^T_0\int_{\errita^d}
p_{T-s}(x-y)\varphi(y)\chi_T(T-s)\\
&&\times\int^s_0{\cal T}_{s-u}
(\varphi\chi(T-u)v_T(\cdot,u))(y)|x|^{-\gamma}dudydsdx\\
I'''_2(T)&=&\frac{H_T}{F^2_T}\int_{\errita^d}\int^T_0\int_{\errita^d}
p_{T-s}(x-y)\varphi(y)\chi_T(T-s)\int^s_0{\cal T}_{s-u}
(v^2_T(\cdot,u))(y)|x|^{-\gamma}dudydsdx.
\end{eqnarray*}
Substituting $u'=u-s, s'=s/T$, and then $x'=xT^{1/\alpha}s^{1/\alpha}$ and 
using part (a) of the Lemma and (\ref{eq:2.24}) we have
\begin{equation}
\label{eq:3.99}
\lim_{T\to\infty}I'_2(T)=\int_{\errita^d}p_1(x)|x|^{-\gamma}dx\int^1_0s^{-\gamma/\alpha}\chi^2(s)ds\int_{\errita^d}\varphi(y)G\varphi(y)dy.
\end{equation}
Applying (\ref{eq:3.5}), (\ref{eq:2.24}) and the Lemma above we get
\begin{eqnarray*}
I''_2(T)&\leq&\frac{CH_T}{F^3_T}\int_{\errita^d}\int^T_0\int_{\errita^d}p_s(x-y)|x|^{-\gamma}\varphi(y)G(\varphi G\varphi)(y)dydsdx\\
&\leq&C_1T^{-(1/2)(1-\gamma/\alpha)}\to 0,
\end{eqnarray*}
and, analogously,
\begin{eqnarray*}
I'''_2(T)&\leq&\frac{C H_T}{F^3_T}\int_{\errita^d}\int^T_0
\int_{\errita^d}p_s(x-y)|x|^{-\gamma}\varphi(y)G((G\varphi)^2)(y)dydsdx\\
&\leq&C_2T^{-(1/2)(1-\gamma/\alpha)}\to 0.
\end{eqnarray*}
This and (\ref{eq:3.98}), (\ref{eq:3.99}) imply that for $\beta=1$ the limit of $({V}/{2})I_1(T)+I_2(T)$ is exactly\\ 
$\log E{\rm exp}\{-C\langle\widetilde{X},\varphi\otimes \psi)\}$. Similar estimations, together with the Lemma,  yield (\ref{eq:3.16}) and (\ref{eq:3.17}).
 This completes the proof of part (a) of the theorem. 
\vglue.5cm
\noindent
{\bf Proof of part (b) of the theorem.} Following the general scheme one can show
\begin{eqnarray*}
\lim_{T\to\infty}I_1(T)&=&\int_{\errita^d}p_1(x)|x|^{-\alpha}dx\chi^{1+\beta}(0)\int_{\errita^d}(G\varphi)^{1+\beta}(y)dy,\\
\lim_{T\to\infty}I_2(T)&=&c_\beta\int_{\errita^d}p_1(x)|x|^{-\alpha}dx\chi^2(0)\int_{\errita^d}\varphi(y)G\varphi(y)dy,
\end{eqnarray*}and 
(\ref{eq:3.16}) and (\ref{eq:3.17})
 (recall that $c_\beta$ is defined by (\ref{eq:2.26})). This is accomplished by an argument similar to the one used in part (a). Due to the criticality $(\gamma=\alpha)$, the integrals $\int^T_0\ldots ds$ in (\ref{eq:3.9})-(\ref{eq:3.11}) require a different treatment. They are split into $\int^1_0\ldots ds+\int^T_1 \ldots ds$; the first summand converges to zero, and in the second one we use the substitution $s'=\log s/\log T$. Here, again, we use repeatedly the Lemma above together with  the easily checked fact that
$$\sup_{T>2}\frac{1}{\log T}\int_{\errita^d}\int^T_0{\cal T}_sh(x)|x|^{-\alpha}dsdx<\infty$$
for any integrable and bounded function $h$ (recall that $d>\alpha)$. We omit details.
\vglue.5cm
\noindent
{\bf Proof of part (c) of the theorem.} Recall that the case $\gamma>d$ has been proved in \cite{BGT6}. For $\alpha<\gamma\leq d$ we use the Lemma (part (c) is particularly important). We show
\begin{eqnarray*}
\lim_{T\to\infty}I_1(T)&=&\int_{\errita^d}G(G\varphi)^{1+\beta}(x)\mu_\gamma(dx)\chi^{1+\beta}(0),\\
\lim_{T\to\infty}I_2(T)&=&c_\beta\int_{\errita^d}G(\varphi G\varphi)(x)\mu_\gamma(dx)\chi^2(0),
\end{eqnarray*}
(\ref{eq:3.16}) and (\ref{eq:3.17}). Here $\mu_\gamma$ is either given by 
(\ref{eq:1.8})
 or, for $\gamma<d$ one can take $\mu_\gamma(dx)=|x|^{-\gamma}dx$. Again, the details are omitted.

\vglue.5cm
\noindent
{\bf 3.8 Proof of Theorem 2.8}
\vglue.5cm

The proof is based on two general lemmas which, hopefully, are of interest by themselves.
\vglue.5cm
\noindent
{\bf Lemma A.} {\it Let $Y$  be an $(d,\alpha,\beta)$-superprocess with $Y_0=\mu$  and $ N$ be the  empirical process of the corresponding branching particle system. If 
for any bounded Borel set $A\subset\erre^d$, 
\begin{equation}
\label{eq:3.100}
P\left[\int^\infty_0Y_t(A)dt<\infty\right]=1,
\end{equation}
then also
\begin{equation}
\label{eq:3.109}
P\left[\int^\infty_0 N_t (A) dt <\infty\right]=1.
\end{equation}
}\vglue.5cm
\noindent
{\bf Proof of Lemma A.}
Let $\zeta$ denote the standard $\alpha$-stable L\'evy process in $\erre^d$, and let $\xi$ be a Markov process with semigroup
$${\cal S}_t\varphi(x)=E_x\left[{\rm exp}\left\{-\int^t_0\psi(\zeta_s)ds\right\}\varphi(\zeta_t)\right],$$
where $\psi$ is a fixed element of $C^\infty_K(\erre^d)$ (bounded support),
$\psi\geq 0$. The process $\xi$ takes values in $\erre^d\cup \{\dagger \}$, where $\dagger$ is a cemetery point where it remains after killing by exp$\{-\int^t_0\psi(\zeta_s)ds\}$. The infinitesimal generator of $\xi$ is
$$A\varphi(x)=(\Delta_\alpha-\psi(x))\varphi(x).$$

Let $Y^\psi$ be a superprocess in $\erre^d$ constructed from $\xi$ and $(1+\beta)$-branching, with $Y^\psi_0=\mu$. The Laplace functional of its occupation time is given by 
\begin{equation}
\label{eq:3.101}
E{\rm exp}\left\{-\int^t_0\langle Y^\psi_s,\varphi\rangle ds\right\}={\rm exp}\{-\langle \mu,u^\psi_\varphi(t)\rangle \},\,\,\varphi\in C^\infty_K(\erre^d),\,\,\varphi\geq 0,
\end{equation}
where $u^\psi_\varphi(x,t)$ is the unique (mild) solution of
\begin{eqnarray}
\label{eq:3.102}
\frac{\partial}{\partial t}u^\psi_\varphi(x,t)&=&(\Delta_\alpha-\psi(x))u^\psi_\varphi(x,t)-\frac{V}{1+\beta}(u^\psi_\varphi(x,t))^{1+\beta}+\varphi(x),\\
u^\psi_\varphi(x,0)&=&0\qquad\qquad\nonumber
\end{eqnarray}
(cf.\ (\ref{eq:1.6})).
The Laplace functional of the occupation time of the process $N$ is given by 
\begin{equation}
\label{eq:3.103}
E{\rm exp}\left\{-\int^t_0\langle N_s,\varphi\rangle ds\right\}={\rm exp}\{-\langle \mu,\;v_\varphi(t)\rangle\},\,\,\, \varphi\in C^\infty_K (\erre^d), \quad \varphi \geq 0,
\end{equation}
where $v_\varphi(x,t)$ is the unique solution of
\begin{eqnarray}
\frac{\partial}{\partial t}v_\varphi(x,t)&=&\Delta_\alpha v_\varphi(x,t)-\frac{V}{1+\beta}(v_\varphi(x,t))^{1+\beta}+\varphi(x)(1-v_\varphi(x,t))\nonumber\\
\label{eq:3.104}
&=&(\Delta_\alpha-\varphi(x))v_\varphi(x,t)-\frac{V}{1+\beta}(v_\varphi(x,t))^{1+\beta}+\varphi(x),\\
v_\varphi(x,0)&=&0\qquad\qquad\nonumber
\end{eqnarray}
(cf.\ (\ref{eq:1.4})).
Equations (\ref{eq:3.102}) and (\ref{eq:3.104}) coincide for $\varphi=\psi$, hence, from (\ref{eq:3.101}), (\ref{eq:3.103}),
\begin{equation}
\label{eq:3.105}
E{\rm exp}\left\{-\int^t_0\langle Y^\psi_s,\psi\rangle ds\right\}=E{\rm exp}\left\{-\int^t_0\langle N_s,\psi\rangle ds\right\}={\rm exp}\{-\langle\mu,u^\psi_\psi(t)\rangle\}.
\end{equation}

The superprocesses $Y\equiv Y^0$ and $Y^\psi$ are obtained as  
(high-density/short-life/small-particle) limits of the same process $N$,  removing first the killed particles in the case of $Y^\psi$. Hence
\begin{equation}
\label{eq:3.106}
E{\rm exp}\left\{-\int^t_0\langle Y^\psi_s,\varphi\rangle ds\right\}\geq E{\rm exp}\left\{-\int^t_0\langle Y^0_s,\varphi\rangle ds\right\}\quad \hbox{\rm for all}\; t,
\end{equation}
for any $\varphi\geq 0$, in particular for $\varphi=\psi$, hence, from 
(\ref{eq:3.105}), (\ref{eq:3.106}), and (\ref{eq:3.101}) with $\psi=0$,
\begin{equation}
\label{eq:3.107}
{\rm exp}\{-\langle \mu,u^\psi_\psi(t)\rangle\}\geq {\rm exp}\{-\langle \mu,u^0_\psi(t)\rangle\}\quad\hbox{\rm for all}\;t.
\end{equation}
Taking $\mu =\delta_x$ in (\ref{eq:3.101}) we see that
$$u^\psi_\psi(x,t)\nearrow u^\psi_\psi(x)\quad{\rm and}\quad u^0_\psi(x,t)\nearrow u^0_\psi(x)\quad{\rm as}\quad t\nearrow\infty,$$
hence, from (\ref{eq:3.107}),
\begin{equation}
\label{eq:3.108}
\langle\mu,u^\psi_\psi\rangle\leq\langle \mu,u^0_\psi\rangle.
\end{equation}

From (\ref{eq:3.100}), (\ref{eq:3.105}), (\ref{eq:3.108}),
\begin{eqnarray*}
1&=&P\left[\int^\infty_0\langle Y_t,\psi\rangle dt<\infty\right]=\lim_{\theta\searrow 0}{\rm exp}\{-\langle\mu,u^0_{\theta\psi}\rangle\}\\
&\leq&\lim_{\theta\searrow 0}{\rm exp}\{-\langle\mu,u^{\theta\psi}_{\theta\psi}\rangle\}=P\left[\int^\infty_0\langle N_t,\psi\rangle dt<\infty\right],
\end{eqnarray*}
so (\ref{eq:3.109}) is satisfied for any bounded set $A\subset \erre^d$ and the lemma is
 proved. $\hfill\Box$

\vglue.5cm
\noindent
{\bf Lemma B.} {\it Let $N$ be the empirical process of the $(d,\alpha,\beta)$-branching particle system with 
locally finite 
initial intensity measure $\mu$. If (\ref{eq:3.109}) is satisfied for any bounded set $A\subset\erre^d$, then 
\begin{equation}
\label{eq:3.110}
P[\sup\{t:N_t(A)>0\}<\infty]=1
\end{equation}
for any bounded set $A$.}

\vglue.5cm
\noindent
{\bf Proof of Lemma B.}
Let $B_R$ be a closed ball in $\erre^d$ with radius $R$ centered at the origin. Let $(t_i,x_i,\tau_i), i=1,2,\ldots,$ be a sequence of random vectors defined as follows for any realization of the branching particle system. First we exclude all the particles which start inside $B_R$ at time $0$ and their progenies. 
Let $t_1$ be the first time any of the remaining particles enters $B_R, x_1$ is the entry point, and $\tau_1$ is the occupation time of the closed ball $B_1(x_1)$ of radius 1 centered at $x_1$ by the tree generated by the entered particle. We exclude this tree from further consideration. Let $t_2$ be the first time after $t_1$ that any of the remaining particles  enters $B_R$, with 
$x_2$ and $\tau_2$ defined analogously as above; and so on. Let $\eta$ denote the total number of first entries $(t_i,x_i,\tau_i), i=1,\ldots,\eta$. We will show that $\eta<\infty$ a.s.. Suppose to the contrary that $P[\eta=\infty]>0$. By construction, $\sum^\eta_{i=1}\tau_i\leq\int^\infty_0 N_t(B_{R+1})dt$, hence $\sum^\eta_{i=1}\tau_i<\infty$ a.s. by (\ref{eq:3.109}). By the strong Markov property and homogeneity of the motion, conditioned on $\{\eta=\infty\}$ the random variables $\tau_i$ are i.i.d.. Hence
$$P\left[\sum\nolimits^\eta_{i=1}\tau_i=\infty|\eta=\infty\right]=1,$$
and this is a contradiction since, as observed above, $P[\sum^\infty_{i=1}\tau_i=\infty]=0$.

Going back to the particles that start inside $B_R$, there are only finitely many of them since $\mu(B_R)<\infty$. 

In conclusion, with probability 1 only finitely many initial particles generate trees that contribute to the occupation time of any given bounded set, and all those trees become extinct a.s. in finite time by criticality of the branching. So (\ref{eq:3.110}) is proved. $\hfill\Box$
\vglue.5cm
Now,
 to prove Theorem 2.8 it suffices to observe that under its assumptions the corresponding superprocess $Y$ suffers local extinction by Theorem $3_\beta$ of \cite{I2}, hence (\ref{eq:3.100}) is clearly satisfied and the theorem follows immediately from the lemmas.

\vglue.5cm
\noindent
{\bf 3.9 Proof of Proposition 2.9}
\vglue.5cm
First observe that it suffices to prove convergence of finite-dimensional distributions. Indeed, in the proof of Theorem 2.1(a) we have shown tightness of $X_T=Z_T-EZ_T$, and the presence of high density was not relevant in that proof. On the other hand, from Proposition 2.1 of \cite{BGT5} it follows easily that the family of deterministic processes $(E\langle Z_T,\varphi\rangle)_{T\geq 1}$ is tight in $C([0,\tau],\erre),\tau>0$. Hence tightness of $Z_T$ follows.

Without loss of generality we assume that $\tau=1$. Fix $0\leq t_1<t_2<\ldots<t_n\leq 1, \varphi_1,\ldots,\varphi_n\in{\cal S}(\erre^d)$, and we may additionally assume that $\varphi_1,\ldots,\varphi_n\geq 0$. In order to show $\Rightarrow_f$ convergence we prove that
\begin{equation}
\label{eq:3.111}
\lim_{T\to\infty}E{\rm exp}\left\{-\sum^n_{k=1}\langle Z_T(t_k),\varphi_k\rangle\right\}={\rm exp}\left\{-\int_{\errita^d}v(x,1)|x|^{-\gamma}dx\right\},
\end{equation}
where $v$ satisfies (2.36) with $\psi$ given by (2.37) for $\theta_k=\int_{\errita^d}\varphi_k(y)dy$ (as explained in \cite{Ta}, the solution of (2.36) is unique.)

For simplicity we consider $\mu_\gamma(dx)=|x|^{-\gamma}dx$ (it will be clear that the limit is the same as for $\mu_\gamma$ given by (1.8)). Also, to simplify the notation we take $\varphi_1=\ldots=\varphi_n=\varphi$. Essentially the same argument can be carried out in the general case.

As in \cite{Ta} and \cite{BGT4} (the possibility to pass from space-time random variable to the present situation) we have 
\begin{equation}
\label{eq:112}
E{\rm exp}\left\{-\sum^n_{k=1}\langle Z_T({t_ k}),\varphi\rangle\right\}={\rm exp}\left\{-\int_{\errita^d}v_T(x,T)|x|^{-\gamma}dx\right\}, 
\end{equation}
where $v$ satisfies (3.4) with $\chi(t)=\sum^n_{k=1}\UNO_{[0,t_k]}(t), 
\chi_T(t)=\chi({t}/{T})$. Formula (3.112) is an analogue of (3.7), and its form is simpler since now we do not subtract the mean.

The right-hand side of (\ref{eq:112}) 
can be written as
$${\rm exp}\left\{-\int_{\errita^d}h_T(x,1)|x|^{-\gamma}dx\right\},$$ 
where
\begin{equation}
\label{eq:3.113}
h_T(x,t)=T^{d/\alpha-\gamma/\alpha}v_T(xT^{1/\alpha},Tt), \,\, 0\leq t\leq 1.
\end{equation}
To prove (3.111) it suffices show that $h_T(\cdot,1)$ converges to $v(\cdot,1)$ in $L^1(\erre^d,|x|^{-\gamma}dx)$; in fact, we will prove that
\begin{equation}
\label{eq:3.114}
h_T\to v\quad{\rm in}\quad C([0,1], L^1(\erre^d,|x|^{-\gamma}dx)).
\end{equation}
By (\ref{eq:3.113}), (\ref{eq:3.4}), (\ref{eq:2.33}) and (\ref{eq:3.20}) we have 

\begin{eqnarray}
h_T(x,t)&=&\int^t_0{\cal T}_{t-s}T^{d/\alpha}\varphi(T^{1/\alpha}\cdot)(x)\chi(1-s)ds\nonumber\\
&-&T^{\gamma/\alpha}\int^t_0{\cal T}_{t-s}(\varphi(T^{1/\alpha}\cdot)h_T(\cdot,s))(x)\chi(1-s)ds\nonumber\\
\label{eq:3.115}
&-&\frac{V}{1+\beta}\int^t_0{\cal T}_{t-s}(h_T(\cdot,s))^{1+\beta}(x)ds.
\end{eqnarray}
In particular, this implies that
\begin{equation}
\label{eq:3.116}
h_T(x,t)\leq CT^{d/\alpha}\int^t_0{\cal T}_s\varphi(T^{1/\alpha}\cdot)ds.
\end{equation}

Let
\begin{eqnarray}
\lefteqn{R_T(x)}\nonumber\\
\label{eq:3.117}
&&\kern -.75cm = \sup_{t\leq 1}\left|\int^t_0{\cal T}_{t-s}(T^{d/\alpha}\varphi(T^{1/\alpha}\cdot))(x)\chi(1-s)ds\right.
\left.-\int^t_0p_{t-s}(x)\chi(1-s)ds\int_{\errita^d}\varphi(y)dy\right|.
\end{eqnarray}
(compare with the first formula on page 851 of \cite{Ta}). We will show that
\begin{equation}
\label{eq:3.118}
\lim_{T\to\infty}\int_{\errita^d}R_T(x)|x|^{-\gamma}dx=0.
\end{equation}

Applying the usual substitutions and the fact that $p_u(x)$ is a decreasing function of $|x|$ we obtain
\begin{eqnarray*}
\lefteqn{\int_{\errita^d}R_T(x)|x|^{-\gamma}dx\leq C\int_{\errita^d}
\int_{\errita^d}\varphi(y)\left|\int^1_0(p_u(x-T^{-1/\alpha})
-p_u(x))du\right||x|^{-\gamma}dydx}\\
&\leq&C\int_{\errita^d}\varphi(y)\int_{\errita^d}|f(x-T^{-1/\alpha}y)-f(x)|dxdy\\&+&C\left(\int_{|x|\leq 1}\int_{\errita^d}\varphi(y)|x|^{-\gamma p}dydx\right)^{1/p}\left(\int_{|x|\leq 1}\int_{\errita^d}\varphi(y)|f(x-T^{-1/\alpha}y)-f(x)|^qdydx\right)^{1/q},
\end{eqnarray*}
where $f$ is defined by (\ref{eq:3.22}), and $p,q>1$ are such that 
$\gamma p<d,(d-\alpha)q<d,\\
1/p+1/q=1$ (such $p$ and $q$ exist since $\gamma<\alpha<d$). 
Hence 
(\ref{eq:3.118}) easily follows from (\ref{eq:3.23}) and (\ref{eq:3.24}).

Next, we will show that
\begin{equation}
\label{eq:3.119}
\lim_{T\to\infty}\int_{\errita^d}\sup_{t\leq 1}\left|T^{\gamma/\alpha}\int^t_0{\cal T}_{t-s}(\varphi(T^{1/\alpha}\cdot)\chi(1-s)h_T(\cdot,s)(x)\right||x|^{-\gamma}dx=0.
\end{equation}
(\ref{eq:3.116}) and (\ref{eq:3.31}) imply that the expression under the lim can be estimated by
\begin{eqnarray*}
T^{\gamma/\alpha}\int_{\errita^d}f_\gamma(y)\varphi(T^{1/\alpha}y)(f*\widetilde{\varphi}_T)(y)dy
&\leq&T^{\gamma/\alpha}||\varphi(T^{1/\alpha}\cdot)||_{\frac{1+\beta}{\beta}}||f||_{1+\beta}||\varphi||_1\\
&\leq&CT^{\gamma/\alpha-(d/\alpha)\beta/(1+\beta)}\to 0
\end{eqnarray*}
(we have used (\ref{eq:2.33}) and $f\in L^{1+\beta}$ by (\ref{eq:3.23}) and
(\ref{eq:3.24})).

To prove (3.114) we check the Cauchy condition, i.e.,
\begin{equation}
\label{eq:3.120}
J(T_1,T_2):=\int_{\errita^d}\sup_{t\leq 1}|h_{T_1}(x,t)-h_{T_2}(x,t)||x|^{-\gamma}dx\to 0\,\,{\rm as}\,\, T_1,T_2\to\infty.
\end{equation}
Using (3.115), (3.118) and (3.119) we have
\begin{equation}
\label{eq:3.121}
J(T_1,T_2)\leq J_1(T_1,T_2)+\frac{V}{1+\beta}J_2(T_1,T_2),
\end{equation}
where $\lim\limits_{T_1,T_2\to\infty}J_1(T_1,T_2)=0$ and
$$
J_2(T_1,T_2)=\int_{\errita^d}\sup_{t\leq 1}\int^t_0{\cal T}_{t-s}
|h^{1+\beta}_{T_1}(\cdot, s)-h^{1+\beta}_{T_2}(\cdot,s)|(x)ds|x|^{-\gamma}dx.$$
By (\ref{eq:3.31}) and the H\"older inequality,
\begin{eqnarray*}
\lefteqn{J_2(T_1,T_2)\leq C\int_{\errita^d}\sup_{t\leq 1}|h^{1+\beta}_{T_1}
(y,t)-h^{1+\beta}_{T_2}(y,t)|dy}\\
&\leq&C_1\left|\left|
\sup_{t\leq 1}|h_{T_1}(\cdot,t)-h_{T_2}(\cdot,t)|\right|\right|_{1+\beta}
\left(\left|\left|\sup_{t\leq 1}|h_{T_1}(\cdot,t)|\right|
\right|^\beta_{1+\beta}
+\left|\left|\sup_{t\leq 1}|h_{T_2}(\cdot,t)|
\right|\right|^\beta_{1+\beta}\right)
\end{eqnarray*}
Using (\ref{eq:3.116}) and the fact that $f\in L^{1+\beta}$ it is easily seen that
$$\sup_{T\geq 1}\left|\left|\sup_{t\leq 1}h_T(\cdot, t)
\right|\right|_{1+\beta}<\infty.$$

To show (\ref{eq:3.120}) it suffices to prove that

\begin{equation}
\label{eq:3.122}
\lim_{T_1,T_2\to\infty}\left|\left|\sup_{t\leq 1}|h_{T_1}
(\cdot,t)-h_{T_2}(\cdot,t)|\right|\right|_{1+\beta}=0.
\end{equation}
This can be derived in a similar way as in \cite{Ta} (see (\ref{eq:2.21}) and subsequent estimates therein). The only difference is that the term corresponding to $I_2(T)$ in \cite{Ta} requires
 a slightly more delicate treatment; in our case it has the form
$$T^{(\gamma/\alpha)(1+\beta)}\int_{\errita^d}
\sup_{s\leq 1}\left|\int^s_0{\cal T}_{s-u}(\varphi(T^{1/\alpha}\cdot)h_{T_1}(\cdot,u))(x)du\right|^{1+\beta}dx.$$
Using (\ref{eq:3.116}) and the H\"older inequality this is estimated by
$$
T^{(d/\alpha)(1+\beta)}||f*(\varphi(T^{1/\alpha}\cdot )(f*\widetilde{\varphi}_T))||^{1+\beta}_{1+\beta}
\leq T^{(\gamma/\alpha)(1+\beta)}||f||^{1+\beta}_{1+\beta}||\varphi(T^{1/\alpha}\cdot)||^{1+\beta}_p||f||^{1+\beta}_q||\varphi||^{1+\beta}_1,
$$
where $q=d/(d-\alpha+\varepsilon), p=d/(\alpha-\varepsilon)$, and $\varepsilon>0$ is such 
that $\gamma<\alpha-\varepsilon$. Then the right-hand side is not bigger than $C T^{((1+\beta)/\alpha)(\gamma+\varepsilon-\alpha)}$,
 which tends to zero as $T\to\infty$.

Combining (\ref{eq:3.120}), (\ref{eq:3.118}), (\ref{eq:3.119}) and 
(\ref{eq:3.122}), it is seen that one can pass to the limit in (\ref{eq:3.115})
letting $T\to\infty$, thus obtaining that the limit of $h_T$ satisfies 
(\ref{eq:2.36}). This proves (\ref{eq:3.114}) and completes the proof of the Proposition. $\hfill\Box$
\vglue.5cm
\noindent
{\bf  Proof of Theorem 2.11}
\vglue.5cm
 The proof is similar to those for the particle system, starting from  an equation analogous to (3.7)-(3.8), where $\widetilde{X}_T$ is now defined for the occupation time fluctuation process (1.10) corresponding to the $(d,\alpha,\beta,\gamma)$ superprocess $Y$, and equation 
(\ref{eq:3.4}) is  replaced by
\begin{equation}
\label{eq:(3.123)}
v_T(x,t)=\int^t_0{\cal T}_{t-u}\left[\varphi_T\chi_T(T-u)-\frac{V}{1+\beta}v^{1+\beta}_T(\cdot,u)\right](x)du,
\end{equation}
where the term $I_2(T)$ in equation (\ref{eq:3.8}), given by
(\ref{eq:3.10}), does not appear. This reflects the fact that comparing the log-Laplace equations (1.4) for the particle system and (1.6) for the superprocess,  the term $-\varphi v_\varphi$ is missing in (1.6).
(Equation (3.123) can be obtained from (3.4) by the same limiting procedure that yields the superprocess from the branching particle system. An equation analogous to (3.7) for the superprocess can be derived from continuous dependence of the occupation time process with respect to the superprocess, and continuity of the mapping $C([0,\tau],{\cal S}'(\erre^d)\ni x\mapsto\widetilde{x}\in {\cal S}'(\erre^{d+1})$ in (1.11) \cite{BGR}.)
 It follows that the results for the superprocess are the same as those for the particle system, except in the cases where $I_2(T)$ has a non-zero limit, and to obtain the results in those cases it suffices to delete those non-zero limits. Therefore the limits in  Theorems 2.1 and 2.5 are the same for the superprocess, and for those in  Theorem 2.6, $c_\beta=0$ in all cases.
$\hfill\Box$
\vglue.5cm
\noindent
{\bf Acknowledgment.} We thank Professors Zenghu Li and Xiaowen Zhou for useful discussions on the superprocess. We are grateful for the hospitality of the Institute of Mathematics, National University of Mexico, where this work was partially done.

\noindent
\def\refname{\hbox{\normalsize\bf References}}

\end{document}